\documentclass[12pt,a4paper,makeidx,reqno]{amsart}

\makeindex \topmargin=0in \headheight=8pt \headsep=0.5in \textheight=8.9in
\usepackage{amsmath,amssymb,amsthm}
\usepackage{mathrsfs}
\usepackage{eucal}
\usepackage[matrix,arrow,curve]{xy}
\usepackage{textcomp}
\usepackage{graphicx}
\usepackage{wrapfig}
\usepackage{verbatim}

\newtheorem{thma}{Theorem}
\newtheorem{lemma}{Lemma}[section]
\newtheorem{cor}{Corollary}
\newtheorem{prop}{Proposition}[section]
\theoremstyle{remark}
\newtheorem{rem}{Remark}[section]
\newtheorem*{defin}{Definition}

\numberwithin{equation}{section}

\def\eps{\varepsilon}

\def\ord{\mathrm{ord}\,}

\def\M{\CMcal{M}_g}
\def\Mp{\CMcal{M}_{g,1}}
\def\Mpp{\CMcal{M}_{g,\{1,1\}}}
\def\S{\CMcal{S}_g}

\def\cM{\overline{\CMcal{M}}_g}

\def\So{\CMcal{S}^-_g}
\def\cSo{\overline{\CMcal{S}}^-_g}

\def\cCo{\overline{\CMcal C}_g}
\def\cSop{\overline{\CMcal{S}}^-_{g,1}}

\def\Sopt{\CMcal{S}^-_{g,2}}

\def\scSo{\overline{\mathcal{S}}^-_g}

\def\tSo{\CMcal{S}^{\circ}_g}
\def\tSop{\CMcal{S}^{\circ}_{g,1}}
\def\tSopt{\CMcal{S}^{\circ}_{g,2}}
\def\tCo{\CMcal{C}_g^{\circ}}
\def\tCop{\CMcal{C}_{g,1}^{\circ}}
\def\Zop{Z_{g,1}}

\def\F{\CMcal{F}}

\def\Ado{\CMcal{H}^-_g}

\def\Pic{\mathrm{Pic}}

\def\dim{\mathrm{dim\,}}
\def\div{\mathrm{div}\,}

\def\Ad{\CMcal{H}_g}

\def\L{\CMcal{L}}

\def\E{\mathbb E_g}
\def\Ep{\mathbb E_{g,1}}

\def\Z{\CMcal{Z}_g}
\def\cZ{\overline {\CMcal{Z}_g}}

\def\S{\CMcal{S}}

\def\cM{\overline{\CMcal{M}}_g}

\def\U{\CMcal{U}}

\def\C{\CMcal{C}}

\def\L{\CMcal{T}}
\def\R{\CMcal{R}}
\def\Q{\CMcal{Q}}

\def\O{\CMcal{O}}

\def\W{\mathrm{W}_g}
\def\Wop{\mathrm{W}_{g,1}}
\def\Wopt{\mathrm{W}_{g,2}}
\def\eWopt{\tilde{\mathrm{W}}_{g,2}}
\def\BW{B_{\mathrm{W}_g}}

\def\Cau{\mathrm{Cau}_g}
\def\BCau{B_{\mathrm{Cau}_g}}
\def\cW{\overline{\mathrm{W}}_g}

\def\cCau{\overline{\mathrm{Cau}}_g}

\def\Caup{\mathrm{Cau}_{g,1}}

\def\Zop{\CMcal{Z}_{g,1}}
\def\Zhat{\widehat{\CMcal{Z}}_{g,1}}

\def\g{\mathfrak{g}}
\def\supp{\mathrm{supp}}

\def\smm{\smallsetminus}
\def\Aut{\mathrm{Aut}}
\def\X{\CMcal{X}}
\def\dhat{\widehat{\Delta}}
\def\pr{\mathrm{pr}}

\def\Efib{\CMcal{E}}
\def\Ffib{\CMcal{F}}
\def\fib{\mathrm{fib}}
\def\Bl{\mathrm{Bl}}

\begin{document}

\title{On the class of caustic on the moduli space of odd spin curves.}
\author{Mikhail Basok }
\date{}
\maketitle

\bigskip

\begin{abstract}
  Let $C$ be a smooth projective curve of genus $g\geq 3$ and let $\eta$ be an odd theta characteristic on it such that $h^0(C,\eta) = 1$. Pick a point $p$ from the support of $\eta$ and consider the one-dimensional linear system $|\eta + p|$. In general this linear system is base-point free and all its ramification points (i.e. ramification points of the corresponding branched cover $C\to\mathbb P^1\simeq \mathbb PH^0(C,\eta+p)$) are simple. We study the locus in the moduli space of odd spin curves where the linear system $|\eta + p|$ fails to have this general behavior. This locus splits into a union of three divisors: the first divisor corresponds to the case when $|\eta+p|$ has a base point, the second one corresponds to theta characteristics which are not reduced at $p$ (and therefore $|\eta + p|$ must have a triple point at $p$) and the third one corresponds to the case when $|\eta + p|$ has a triple point different from $p$. The second divisor was studied by G. Farkas and A. Verra in~\cite{FARo} where its expansion in the rational Picard group was used to prove that the moduli space of odd spin curves is of general type for genus at least $12$. We call the first divisor a Base Point divisor and the third one a Caustic divisor (following Arnold terminology for Hurwitz spases). The objective of this paper is to expand these two divisors via the set of standard generators in the rational Picard group of the moduli space of odd spin curves.
\end{abstract}

\section{Introduction.}

In this paper we study several divisor classes in the moduli space $\cSo$ of odd spin curves. This moduli space presents a rare example when the birational type is known for all $g$: a complete classification was established by G. Farkas and A. Verra~\cite{FARo}, where they showed that $\cSo$ is uniruled when $g\leq 11$ and of general type for other values of $g$. As usual, to prove that a moduli space is of general type one needs to show that the canonical class lies in the interior of the cone of effective divisors, or, equivalently, to find a divisor with the slope strictly less than the one of the canonical class. In the case of $\cSo$ the authors of~\cite{FARo} show that a linear combination of the Brill-Noether divisor and a divisor $\cZ$ defined below has the slope small enough. The divisor $\cZ$ is defined as the closure in $\cSo$ of the locus consisting of spin curves $(C, \eta)$ with non-reduced $\eta$, i.e.
\begin{equation*}
  \cZ = \text{closure of }\left\{ (C, \eta)\in \So\ \mid\ \eta = \O_C(2x_1+ x_2+ \dots + x_{g-2}) \right\}.
\end{equation*}
The following relation holds for the class of this divisor in $\Pic(\cSo)\otimes\mathbb Q$ (see~\cite[Theorem~0.5]{FARo}):

\begin{thma}
  One has
  \begin{equation*}
    [\cZ] = (g+8)\,\lambda - \frac{g+2}{4}\alpha_0 - 2\beta_0 - \sum_{i = 1}^{g-1} 2(g-i)\,\alpha_i,
  \end{equation*}
  in $\Pic(\cSo)\otimes \mathbb Q$ where $\lambda,\beta_0, \alpha_0,\alpha_1,\dots, \alpha_{g-1}$ are the standard generators of $\Pic(\cSo)$.
  \label{thma:Farkas_formula}
\end{thma}
\noindent See Section~\ref{sec:moduli_of_spin_curves} for precise definitions.

As we will see, the divisor $\Z$ arises naturally as a component of a pullback of degeneration divisor from the Hurwitz space. Namely, let $(C,\eta)\in \So$ be a generic point and assume that $\eta\simeq \O_C(p+x_1+\dots+x_{g-2})$ (let us write $\eta\geq p$ in this case). As $(C,\eta)$ is chosen generic and $K_C-\eta$ is equal to $\eta$ in the Picard group of $C$ we have $h^0(C,\eta+p) = 1 + h^0(C,\eta-p) = 2$ by Riemann-Roch. Using the general choice of $(C,\eta)$ again we can assume that the linear system $|\eta+p|$ has no base points (see~Lemma~\ref{lemma:W_is_a_divisor} where we show that this condition holds generically). It follows that the linear system $|\eta+p|$ gives a rise to a branched cover of $\mathbb P^1$ that is explicitly defined as 
\begin{equation}
  \begin{split}
    C\to \mathbb PH^0(C,\eta+p)\simeq \mathbb P^1\\
    q\mapsto \{v\in H^0(C,\eta+p):\ v(q) = 0\}.
  \end{split}
  \label{eq:introduction_morphism_to_P1}
\end{equation}
Moreover, for a generic $(C,\eta)$ the linear system $|\eta+p|$ behaves like a generic point from the space $\g_g^1$ parametrizing one-dimensional linear systems on genus $g$ curves, i.e. all ramification points of the corresponding branched cover are simple (i.e. have the ramification index equal to $2$). It means that if a spin curve $(C,\eta)$ belongs to a certain open subset of $\So$ then the branched cover~\eqref{eq:introduction_morphism_to_P1} represents a point in the main stratum of the Hurwitz space of degree $g$ covers of $\mathbb P^1$.

Let us now consider the locus in $\So$ where $|\eta + p|$ fails to have this general behavior for some $p\leq \eta$. This locus splits into three components. The first two correspond to the case when the order one of some ramification point of $|\eta+p|$ increases. Notice that $p$ is in particular a ramification point of $|\eta+p|$. Thus if $(C,\eta)\in \Z$ and $\eta$ is not reduced at $p$ then $|\eta+p|$ has a ramification of the order at least $3$ at $p$, therefore $\Z$ is contained in the locus that we consider. Another possibility is that $|\eta+p|$ has a ramification point of order three that does not belong to $\supp(\eta)$.

\begin{defin}
  Let $g\geq 3$ be chosen. The \emph{Caustic divisor} $\Cau$ is the closure of the locus in the moduli space $\So$ parametrizing odd spin curves $(C,\eta)$ such that there exist $p\in \supp(\eta)$ and $q\in C-\{p\}$ such that $h^0(C,\eta + p- 3q) >0$.
\end{defin}
In the case when $h^0(C,\eta) = 1$ the condition $h^0(C,\eta+p-3q)>0$ is equivalent to the fact that the mapping~\eqref{eq:introduction_morphism_to_P1} has a ramification point of order at least $3$. It is not hard to check that this locus in indeed a divisor, i.e. it has codimension $1$ if non-empty. Finally, it may happen that $|\eta+p|$ has a base point (for example, if $g\geq 3$, $C$ is hyperelliptic, and $h^0(C,\eta) = 1$). The corresponding locus will be denoted by 
\begin{equation*}
  \W = \left\{ (C, \eta)\in \So\ \mid\ \text{$|\eta + p|$ has a base point} \right\}.
\end{equation*}
It is straightforward to show that $\W$ is a divisor if it does not coincide with $\So$. We will call $\W$ the {\it Base Point divisor}. Let $\cCau$ and $\cW$ be closures of $\Cau$ and $\W$ in $\cSo$. The main result of the paper is the following: 

\begin{thma}
  For any $g\geq 4$, $\cCau$ and $\cW$ are divisors in $\cSo$ and their classes in $\Pic(\cSo)\otimes \mathbb Q$ satisfy relations:
  \begin{equation}
    [\cCau] = \frac{9g^2 + 179g - 134}2 \lambda - \frac{9g^2 + 59g - 50}8 \alpha_0 - (24g - 22)\beta_0 - \sum_{j = 1}^{g-1} (g-j)(9g + 27j - 19)\alpha_j,
    \label{eq:caustic_divisor_formula}
  \end{equation}
  \begin{equation}
    [\cW] = \frac 12 \left (  \frac{g^2 + 11g - 6}2 \lambda - \frac{g^2 + 3g - 2}8 \alpha_0 - (2g-2) \beta_0 - \sum_{j = 1}^{g-1} (g - j)(g + 3j - 3)  \alpha_j  \right),
    \label{eq:base_point_divisor_formula}
  \end{equation}
  where $\lambda, \beta_0, \alpha_0,\alpha_1,\dots, \alpha_{g-1}$ are the standard generators of $\Pic(\cSo)\otimes \mathbb Q$.
  \label{thm:caustic_and_base_point_divisors_formula}
\end{thma}

\begin{proof}
  Denote the coefficients in expansions of $[\cW]$ and $[\cCau]$ by the formulas:
  \begin{align*}
    & [\cCau] = l^c \cdot \lambda - a_0^c \cdot \alpha_0 - b_0^c\cdot \beta_0 - \sum_{j = 1}^{g-1}a_j^c \cdot \alpha_j,\\
    & [\cW] = l^w \cdot \lambda - a_0^w \cdot \alpha_0 - b_0^w\cdot \beta_0  - \sum_{j = 1}^{g-1}a_j^w \cdot \alpha_j.
  \end{align*}
  By Proposition~\ref{prop:relation_between_W_and_lambda} and Proposition~\ref{prop:relation_between_Cau_and_lambda} we have
  \begin{equation*}
    \begin{split}
      & l^c =  \frac{9g^2 + 179g -134}{2},\\
      & l^w = \frac{g^2 + 11g - 6}{4}.
    \end{split}
  \end{equation*}
  By Proposition~\ref{prop:alphaj_for_j>0}, for each $j = 1,\dots, g-1$ we have:
  \begin{align*}
    &a_j^c = (g-j)(9g + 27j - 19)\\
    &a_j^w = \frac{1}{2}\cdot(g-j)(g + 3j - 3).
  \end{align*}
  Finally, Corollary~\ref{cor:alpha0_and_beta0} implies that
  \begin{align*}
    & a_0^w = \frac{g^2 + 3g - 2}{8},\qquad b_0^w = g-1,\\
    & a_0^c = \frac{9g^2 + 59g - 50}{8},\qquad b_0^c = 24g-22.
  \end{align*}
\end{proof}

\noindent The precise definition of $\lambda, \beta_0, \alpha_0,\alpha_1,\dots, \alpha_{g-1}$ are given in Section~\ref{sec:moduli_of_spin_curves}. Note that one can alternatively define $\W$ by
\begin{equation*}
  \W = \left\{ (C, \eta)\in \So\ \mid\ \text{there exist } p,q\in \supp(\eta) \text{ s.t. } h^0(\eta + p- q) > 1 \right\}.
\end{equation*}
The condition $h^0(C, \eta + p - q) > 1$ is symmetric in $p$ and $q$, which explains the global factor~$\frac 12$ in the right-hand side of~\eqref{eq:base_point_divisor_formula}.

Now let us say a few words about the strategy of the proof of Theorem~\ref{thm:caustic_and_base_point_divisors_formula}. Computation of the coefficients splits into two main parts. In the first part we establish a relation between $\lambda$ and the both divisors $[\cW]$ and $[\cCau]$ working over the open part $\So$ of $\cSo$. The main observation here is the following. Let $\pi:\C\to B, \L\to\C$ be a smooth family of odd spin curves (see~Section~\ref{subsec:universal_family} for the definition of a family of spin curves). For $b\in B$ denote by $(C_b,\eta_b)$ the spin curve $(\pi^{-1}(b), \L|_{\pi^{-1}(b)})$. If $p:B\to \C$ is a section of $\pi$ such that $\eta_b\geq p_b$ for each $b$, consider the variety
\begin{equation}
  V = \{(b,q_b)\ \mid\ q_b\in C_b:\ \eta_b + p_b \geq 2q_b,\ q_b\neq p_b\}.
\end{equation}
Then the forgetful map $V\to B$ has the degree $4g-3$ and its ramification locus is the union $B_{\W}\cup B_{\Cau}$, where
\begin{equation*}
  \begin{split}
    & B_{\W} = \{ (b,q_b)\ \mid\ q_b\text{ is a base point of }|\eta_b + p_b| \}\\
    & B_{\Cau} = \{(b,q_b)\ \mid\ \eta_b+p_b\geq 3q_b,\ q_b\neq p_b\}.
  \end{split}
\end{equation*}
Then we have
\begin{equation}
  [B_{\W}] + [B_{\Cau}] = (\omega_{\pi} + [V])\cdot [V]
  \label{eq:introduction_adjuction_formula}
\end{equation}
where $\omega_{\pi}$ is the relative dualizing sheaf (this formula assumes that $V$ is smooth and the ramification is simple; these technical issues will be carefully resolved in Section~\ref{sec:construction_over_the_smooth_part}). To compute the right-hand side of~\eqref{eq:introduction_adjuction_formula} we use the fact that $V$ is a component of the ramification locus of a certain rational morphism $\C\dasharrow P$, where $P\to B$ is a $\mathbb P^1$-bundle with fiber $\mathbb P H^0(C,\eta_b+p_b)$ over a point $b\in B$. The morphism is defined fiber-wise by the rule~\eqref{eq:introduction_morphism_to_P1}. Using these construction we relate classes $\pi_*[C_{\W}]$ and $\pi_*[B_{\Cau}]$ with several other classes associated with $B$ and $p$. As the relation we obtain functorially depends of the family of spin curves we actually get a relation in $\Pic(\So)\otimes \mathbb Q$ between $[\Cau]$, $[\W]$ and $\lambda$. Another relation between these classes follows from relating the blow up of $\C$ along $B_{\W}$ with the variety
\begin{equation*}
  P_1 = \{(b, q_b, \sigma)\ \mid\ \sigma\in \mathbb PH^0(C,\eta_b+p_b),\ \sigma(q_b) = 0\}.
\end{equation*}
The variety $P_1$ arises as a divisor in the fibered product $P\times_B \C$ and the projection $P_1\to \C$ is isomorphic to the blow up of $\C$ along $B_{\W}$. Comparing $K_{P_1}$ with $K_{\C}$ and using the adjunction formula for $P_1$ we relate $[\W]$ with $\lambda$.

The second part of the proof of the main theorem deals with computation of all the boundary coefficients. Here we use the classical machinery of the intersection theory, namely, we compute intersection numbers of $\cW$ and $\cCau$ with several test families and obtain sufficient number of equations to reconstruct all the coefficients.

The computations that we briefly outlined above turn out to be quite lengthy after carefully writing all the details. The resulting formulas in Theorem~\ref{thm:caustic_and_base_point_divisors_formula} are not looking very simple either. It is therefore natural to look for a test to check the result. We find a relation between $2[\W] + [\Cau]$ and $\lambda$ using an alternative approach in Section~\ref{sec:verification_section}. The computation of intersection numbers of $\cW$ and $\cCau$ with certain test families also provides a verification of our formulas (see Remark~\ref{rem:test_via_intersection}).

\bigskip
The paper is organized as follows. In Section~\ref{sec:moduli_of_spin_curves} we recall some basic facts about the moduli space of odd spin curves and its compactification constructed by Cornalba. In Section~\ref{sec:construction_over_the_smooth_part} we consider another spaces that parametrize spin curves with several additional data. As a result we end up with relations between $\lambda$ and $[\W]$ and $[\Cau]$. In Section~\ref{sec:test_curves} we compute intersection numbers of $\cW$ and $\cCau$ with several test families and this completes the preparation for the proof of Theorem~\ref{thm:caustic_and_base_point_divisors_formula}. The computations of Section~\ref{sec:test_curves} also provide a verification of the boundary contribution in the formulas in Theorem~\ref{thm:caustic_and_base_point_divisors_formula}. Finally, in Section~\ref{sec:verification_section} we verify the relation between $\lambda$ and $[\W]$ and $[\Cau]$.

\subsection*{Acknowledgements} I would like to thank Gavril Farkas for stating the problem and for discussion of my work in Saint-Petersburg. I am also very grateful to Alexander Kuznetsov for several discussions that helped me to understand important details. Finally, I would like to thank my advisor Peter Zograf.

\section{Moduli space of odd spin curves}\label{sec:moduli_of_spin_curves}

Let $\M$, $g\geq 3$, be the moduli space of smooth genus $g$ algebraic curves. Let $\cM$ be its Deligne-Mumford compactification. The boundary $\cM\smallsetminus\M$ consists of $\left [ \frac g2 \right ] + 1$ irreducible divisors $\Delta_0,\dots, \Delta_{\left [ \frac g2 \right ] }$ where $\Delta_0$ is the closure of the locus of irreducible curves with one node and $\Delta_j$ for $j\geq 1$ is the closure of the locus of reducible one-nodal curves such that two irreducible components are of genuses $j$ and $g-j$.

The moduli space $\So$ of {\it smooth} odd spin curves is a degree $2^{g-1}(2^g-1)$ cover of $\M$. The cover is extended to a branched cover of $\cM$ by the Cornalba compactification $\cSo$ of $\So$ ramified over~$\Delta_0$. Let us briefly discuss some basic facts about this compactification (see~\cite{COR} for all details).

\subsection{Cornalba compactification.} Given a nodal curve $C$ we call a rational component $E\subset C$ of it {\it exceptional} if $E\cap \overline{C \smallsetminus E_i} = 2$. A nodal curve $C$ is called {\it quasi-stable} if it satisfies two conditions:

1) Any rational component $E$ of $C$ intersects $\overline{C \smallsetminus E}$ at two or more points;

2) Any two exceptional components of $C$ are disjoint.

Following~\cite{COR} we define a {\it spin curve} as a triple $(C,\eta,\beta)$ consisting of a quasi-stable curve $C$, a line bundle $\eta$ of degree $g-1$ on it and a homomorphism $\beta:\eta^{\otimes 2}\to \omega_C$ with the following properties:

1) $\eta$ is of degree one on every exceptional component of $C$;

2) $\beta$ is not zero on every non-exceptional component of $C$.

\noindent It follows from this properties that $\eta^{\otimes 2}|_{C_1}$ is isomorphic to $\omega_{C_1}$, where $C'\subset C$ is the union of all non-exceptional components of $C$. The parity of the spin curve $(C,\eta,\beta)$ is equal by definition to the parity of $\dim H^0(C,\eta)$. The parity is invariant under continuous deformations (see~\cite{MUM} and~\cite{ATJ}). An isomorphism between $(C,\eta,\beta)$ and $(C',\eta',\beta')$ is an isomorphism $\sigma~:~C\to~C'$ and an isomorphism $\phi$ between $\sigma^*\eta'$ and $\eta$ such that the following diagram
$$
\xymatrix{
\eta^2 \ar[rr]^(0.45){\phi\otimes\phi} \ar[d]^{\beta}&& \sigma^*(\eta' )^2\ar[d]^{\sigma^*\beta'}\\
\omega_C\ar[rr]^(0.45){\simeq} && \sigma^* \omega_{C'}
}
$$
is commutative. In~\cite{COR} Cornalba construct an analytic structure on the equivalence classes of all spin curves of given parity and genus. We denote by $\cSo$ the corresponding moduli space of odd spin curves of genus $g$. There is a natural map $\rho~:~\cSo~\to~\cM$ which maps (an equivalence class of) a triple $(C,\eta,\beta)$ to (an equivalence class of) a curve $C'$, where $C'$ is obtained from $C$ by contracting all exceptional components to points.

\subsection{Rational Picard group of $\cSo$.} Let $\rho: \cSo\to \cM$ be the forgetful map. The boundary $\cSo\smallsetminus\So$ is the union of irreducible divisors $A_0,\dots,A_{g-1}, B_0$ such that $\rho(A_0) = \rho(B_0) = \Delta_0$ and $\rho(A_{g-j}) = \rho (A_j) = \Delta_j$ for $j = 1,\dots,\left [ \frac g2 \right ]$.

{\it Description of $A_j$ for $j \neq 0$.} Note that there are no spin curves $(C,\eta,\beta)$ with a reducible one-nodal base curve $C$, since the relative dualizing sheaf $\omega_C$ on a reducible curve with one node being restricted to each component must be of odd degree (see~\cite{COR},~\cite[p.5]{FARo} for more details).

Let $(C,\eta,\beta)$ be a spin curve such that $C = C_-\cup E\cup C_+$ where $C_-$ and $C_+$ are smooth curves of genus $j$ and $g-j$ respectively and $E$ is an exceptional component. The divisor $A_j$ parametrizes the closure of the locus of such curves with the property that $\eta$ restricted to $C_-$ is odd.

{\it Description of $A_0$ and $B_0$.} Contrary to the case $j\neq 0$, a spin curve $(C,\eta,\beta)$ for which $\rho(C,\eta,\beta)$ is an irreducible one-nodal curve, does not necessarily have exceptional components. Let $A_0$ parametrize the closure of the locus of spin curves with one-nodal irreducible underlying curve and $B_0$ parametrize the closure of the locus of reducible spin curves (with an exceptional component) that are mapped to $\Delta_0$ by $\rho$.  It can be checked that the map $\rho$ has simple ramification along $B_0$ and is unramified on $\cSo\smallsetminus B_0$.

Denote by $\alpha_j$ the class of $A_j$ and by $\beta_0$ the class of $B_0$ in the rational Picard group $\Pic(\cSo)\otimes \mathbb Q$ respectively. Let $\lambda$ be the pullback of the Hodge class on $\cM$ under $\rho$. The rational Picard group is generated by these classes:

\begin{equation*}
  \Pic(\cSo)\otimes \mathbb Q = \mathrm{span}_{\mathbb Q}(\lambda, \beta_0, \alpha_0,\alpha_1,\dots, \alpha_{g-1}).
\end{equation*}

\section{Geometric interpretation of $\W$ and $\Cau$.}\label{sec:construction_over_the_smooth_part}

Starting from this point we assume that $g\geq 4$. For such values of $g$ the open subvariety $\tSo$ of $\So$ consisting of $(C,\eta)$ such that
\begin{equation}
  \begin{split}
    & h^0(C, \eta) = 1,\\
    & |\Aut(C)| = 1
  \end{split}
  \label{eq:no_automorphisms_and_h0=1}
\end{equation}
is a complement to a variety of codimension at least $2$. The condition $|\Aut(C)|=1$ ensures that $\tSo$ is smooth. On the other hand we have $\Pic(\tSo) = \Pic(\So)$, thus we can study $\tSo$ instead of $\So$ in order to relate $\lambda$ and the classes of $\Cau$ and $\W$. Recall that given a curve $C$ and an isomorphism class of a line bundle $L\in \Pic(C)$ and an effective divisor $D$ on $C$ we write $L\geq D$ if there exists a divisor $D'$ on $C$ such that $L\simeq \O_C(D')$ and $D'\geq D$.

Our goal is to construct a variety parametrizing the set
\begin{equation}
  \left\{ (C, \eta)\in \tSo, p, q\in C\ \mid\ \eta\geq p,\ \eta+p\geq 2q \right\}.
  \label{eq:definition_of_Sopt_as_a_set}
\end{equation}
This variety is expected to have two irreducible components: the first one corresponds to $p = q$ and the second component is of our particular interest. Denote it by $\tSopt$ in advance. The direct image of the ramification divisor of the forgetful map $\tSopt\to \So$ is a linear combination of classes $[\W]$, $[\Cau]$. Our strategy is to use this observation to express the classes of $\Cau$ and $\W$ via $\lambda$.

\subsection{Moduli space of spinors with marked zero.}
\label{sec:spinors_with_marked_zero}

As a first step to our goal we construct and study a variety parametrizing the set
\begin{equation*}
  \{(C,\eta, p)\ \mid\ \eta\geq p,\ p\in C \}.
\end{equation*}
Let $l:\tCo\to \tSo$ be the universal spin curve over $\tSo$. We say that an invertible sheaf $\F$ on $\tSo$ is a \emph{universal spin line bundle} if for any $(C,\eta)\in \tSo$ the restriction of $\F$ to $l^{-1}([C,\eta])$ is isomorphic to $\eta$. It is straightforward to show that such an $\F$ exists, but the properties stated above do not define $\F$ uniquely: indeed, if $\F_0$ is any invertible sheaf on $\tSo$ then $l^*\F_0\otimes \F$ is again a universal spin line bundle. Note that, conversely, if $\F_1$ and $\F_2$ are two universal spin line bundles, then $\F_1\otimes \F_2^{-1}$ is isomorphic to a pullback of a sheaf from $\tSo$. Thus if we additionally require that $c_1(l_*(\F)) = 0$ then the universal spin line bundle will be defined uniquely up to an isomorphism. Note that $l_*\F$ is a (trivial) line bundle with fibers isomorphic to $H^0(C,\eta)$. Moreover, since $\omega_l\otimes \F^{-1}$ is again a universal spin line bundle then $\omega_l\otimes \F^{-1}\simeq \F\otimes l^*\F_0$ for some $\F_0$ on $\tSo$. It follows that $R^1l_* \F \simeq \F_0^{-1}$. Using this observation and Grothendieck-Riemann-Roch formula we find that $c_1(\F_0) = -\frac{\lambda}{2}$ and thus
\begin{equation}
  c_1(\F) = \frac{1}{2}c_1(\omega_l)+\frac{\lambda}{4}.
  \label{eq:first_chern_class_of_spin_bundle}
\end{equation}
Let $\tSop\subset \tCo$ be the zero locus of the line bundles homomorphism $l^*l_* (\F)\to \F$. A generic $(C,\eta)\in \tSo$ has the property that $\eta\simeq \O_C(D)$ where $D$ an effective divisor on $C$ such that $\ord_x D\leq 1$ for each $x\in C$, therefore $\tSop = \div(l^*l_* (\F)\to \F)$, in particular $\tSop = c_1(\F)$. Let
\begin{equation}
  \mu := l_*\left( c_1(\omega_l)\cdot \tSop \right).
  \label{eq:definition_of_mu}
\end{equation}

\begin{lemma}
  We have
  \begin{equation*}
    \mu = \frac{g+11}{2}\, \lambda
  \end{equation*}
  in $\Pic(\tSop)\otimes \mathbb Q$.
  \label{lemma:relation_between_mu_and_lambda}
\end{lemma}

\begin{proof}
  Using that $\tSop \equiv c_1(\F)$ (where $\equiv$ means the equality of the class of $\tSop$ in $\Pic(\tCo)\otimes \mathbb Q$) and~\eqref{eq:first_chern_class_of_spin_bundle} we can write
  \begin{equation*}
    \mu \equiv l_*\left( c_1(\omega_l)\cdot c_1(\F) \right) \equiv  l_*\left( c_1(\omega_l)\cdot \left(\frac{1}{2}c_1(\omega_l)+\frac{\lambda}{4}\right) \right)
  \end{equation*}
  Using the Mumford formula that reads as $l_*(\omega_l^2) \equiv 12\lambda$ in our situation we write
  \begin{equation*}
    \mu \equiv l_*\left( c_1(\omega_l)\cdot \left(\frac{1}{2}c_1(\omega_l) + \frac{\lambda}{4}\right) \right) \equiv \frac{g+11}{2}\, \lambda.
  \end{equation*}
\end{proof}

Let us prove that $\tSop$ is smooth outside a subvariety of codimension $2$. In this purpose we need to recall the notion of Scorza curve. Recall that given a smooth non-rational curve $D$ and an even theta characteristic $\eta_D$ on $D$ the associated Scorza curve is defined as
\begin{equation}
  T_{\eta_D} := \{(p,q)\in D\times D\ \mid\ \eta_D + p \geq q\}
  \label{eq:def_of_Scorza}
\end{equation}

\begin{lemma}
  Let $D$ be a generic curve of genus $g-1$ and $\eta_D$ be an even theta characteristic on $D$. Then there exists an injective morphism $\iota: T_{\eta_D} \to \cSop$, where $\cSop$ is the closure of $\tSop$ in the universal spin curve $\CMcal{C}\to \cSo$. Moreover, the image of $\iota$ is smooth in $\CMcal C$.
  \label{lm:Scorza_into_tSop}
\end{lemma}

\begin{proof}
  By~\cite[Theorem~4.1.]{FARo} the smoothness of $T_{\eta_D}$ is a generic property, so we can choose $D$ such that the curve $T_{\eta_D}$ is smooth. Choose an arbitrary smooth elliptic curve $X$ with a marked point $x\in X$ and consider the family
  \begin{equation*}
    G_1 = \left\{ X\cup_x E\cup_p D, q\in D, \eta_X, \eta_D \right\}_{p\in D}\subset \CMcal C.
  \end{equation*}
  Define $\iota : T_{\eta_D}\to G_1$ by $(p,q)\mapsto  (X\cup_x E\cup_p D, q\in D, \eta_X, \eta_D)\in G_1$. We have $\iota(T_{\eta_D}) = G_1\cap \overline{\tSop}$. A straightforward check shows that $\iota$ satisfies all desired properties. 
\end{proof}

Note that the forgetful morphism $\tSop\to \tSo$ is a branched cover. Denote by $\Zop$ the ramification locus of it; alternatively $\Zop$ can be defined as
\begin{equation}
  \Zop := \left\{ (C, \eta, p)\in \tSop\ \mid\ \eta\geq 2p \right\}.
  \label{eq:definition_of_Zop}
\end{equation}

\begin{lemma}
  The variety $\tSop$ is smooth outside a subvariety of codimension at least~$2$.
  \label{lm:tSop_is_smooth}
\end{lemma}

\begin{proof}
  Obviously, the singular locus of $\tSop$ is contained in $\Zop$. Our first step is to show that $\Zop$ is irreducible. We deduce it from a famous result~\cite{ConComp} on connectedness of the moduli space of translation surfaces. Namely, let $\E\to \M$ be the Hodge bundle and let $\mu = (\mu_1,\dots,\mu_n)$ be a collection of positive integers such that $\sum_{i = 1}^n\mu_i = 2g-2$. Let $\Ad(\mu)\subset \mathbb P\E$ denote the open variety given by
  \begin{equation*}
    \Ad(\mu) = \{ (C,\omega)\ \mid\ \exists x_1,\dots, x_n\in C:\  \div \omega = \mu_1x_1 + \mu_2 x_2 +\dots + \mu_nx_n,\ x_i\neq x_j \text{ if } i\neq j \}.
  \end{equation*}
  Then the variety $\Ad(\mu)$ is a smooth orbifold for each $\mu$ and connected components of $\Ad(\mu)$ are completely classified in~\cite{ConComp}. For our purpose we will consider $\mu = \{4,2,2,\dots,2\}$ where we denote $\{4,2,2,\dots,2\}$ by $\{4, 2^{g-3}\}$ for simplicity. Any differential $\omega$ on a curve $C$ that represents a point $[C,\omega]$ from the space $\Ad(4,2^{g-3})$ defines a spin structure $\eta \simeq \O_C(\frac{1}{2}\div \omega)$. The parity of the spin structure is a topological invariant, so let $\Ado(4,2^{g-3})$ denote the subvariety of $\Ad(4,2^{g-3})$ consisting of those $\omega$ that define an odd spin structure. Then results of~\cite{ConComp} imply that the variety $\Ado(4,2^{g-3})$ is smooth and connected orbifold for $g\geq 4$.
  
  Now, given a point $(C,\eta, p)\in \Zop$ we consider the holomorphic differential $\omega_{\eta}$ on $C$ that corresponds to the square of a unique (up to a multiplicative constant) non-zero element in $H^0(C,\eta)$ under the isomorphism $\eta^{\otimes 2}\simeq K_C$. Consider the map $\Zop\to [\text{closure of } \Ado(4,2^{g-3})]$ that sends $(C,\eta,p)$ to $[C,\omega_{\eta}]$. This map has the degree $1$, since a generic $\omega_{\eta}$ has a unique zero of order $4$. Therefore the connectedness of $\Ado(4,2^{g-3})$ implies that $\Zop$ is irreducible.
  
  Using the irreducibility of $\Zop$ we conclude that the singular locus of $\tSop$ either coincides with $\Zop$ or has codimension bigger than $1$. Let us show that the first case cannot occur. For this purpose we will use the map $\iota$ constructed in Lemma~\ref{lm:Scorza_into_tSop}. Suppose that $(\tSop)^{\mathrm{sing}} = \Zop$ and consider the closure $\overline{\Zop}$ of $\Zop$ in $\CMcal C$. Then it is easy to check that the range of $\iota$ intersects $\overline{\Zop}$. But then the range of $\iota$ cannot be smooth, since it is contained in the closure of $\tSop$. Thus we obtained a contradiction.

  It remains to show that $\tSop$ is irreducible. Here we again use the map $\iota$. Suppose that $\tSop$ has two irreducible components, say, $\CMcal{S}_1$ and $\CMcal{S}_2$. Note that both maps $\CMcal{S}_1, \CMcal{S}_2\to \tSo$ must be surjective, otherwise the map $\tSop\to \tSo$ would have degree bigger then $g-1$ at some point. Thus closures $\overline{\CMcal{S}}_1$ and $\overline{\CMcal{S}}_2$ intersect $\iota(T_{\eta})$. Since the degree of the map $\tSop\to \tSo$ is equal to the degree of the map $\iota(T_{\eta})\to \cSo$ we see that $\overline{\CMcal{S}}_1\cap \iota(T_{\eta})$ and $\overline{\CMcal{S}}_2\cap \iota(T_{\eta})$ cannot coincide, so we conclude that the Scorza curve $T_{\eta}$ is reducible. But this cannot be the case since a generic Scorza curve is smooth and connected (for the connectedness of Scorza curves see~\cite[Section~7.1]{DK}).
\end{proof}

\subsection{Muduli space of spinors having a pole}

Let $\tCop = \tSop\times_{\tSo}\tCo$ be the pullback of the universal spin curve to $\tSop$, set theoretically 
\begin{equation}
  \tCop = \{(C,\eta,p,q)\ \mid\ (C,\eta,p)\in \tSop,\quad q\in C\}.
  \label{eq:DefOftCop}
\end{equation}
Let $\pi: \tCop\to\tSop$ be the forgetful morphism. By abusing of notations we denote the pullback of $\F$ to $\tCop$ by $\F$. Recall that the choice of $\F$ ensures that $F := \pi_*\F$ is isomorphic to $\O_{\tSop}$ and by~\eqref{eq:first_chern_class_of_spin_bundle} we have
\begin{equation}
  c_1(\F) = \frac{1}{2}c_1(\omega_{\pi})+\frac{\lambda}{4}.
  \label{eq:c_1_F}
\end{equation}
The morphism $\pi$ has a diagonal section $\delta: \tSop\to\tCop$ that sends $(C, \eta, p)$ to $(C, \eta, p, p)$. We denote its image by $\Delta$ and consider the sheaf 
\begin{equation}
  E = \pi_*\F(\Delta).
  \label{eq:DefOfE}
\end{equation}
The sheaf $E$ is a locally free sheaf of rank $2$ with the fiber over $(C,\eta,p)$ isomorphic to $H^0(C,\eta\otimes \O_C(p))$. Note that elements of $H^0(C,\eta\otimes \O_C(p))$ can be thought of as the space of sections of $\eta$ having simple pole at $p$. Applying $\pi_*$ to the short exact sequence
\begin{equation*}
  0\to \F\to \F(\Delta) \to \F(\Delta)|_{\Delta} \to 0
\end{equation*}
we obtain the long exact sequence
\begin{equation*}
  0\to \O_{\tSop}\to E\to \pi_*(\F\otimes\omega_{\pi}^{\vee}|_{\Delta}) \to R^1\pi_*\F \to  R^1\pi_*\F(\Delta) \to 0.
\end{equation*}
Since $\pi_*\F$ and $\pi_*\F(-\Delta)$ are line bundles and $c_1(\omega_{\pi}\otimes \F^{\vee}) = c_1(\F) -\frac{\lambda}{2}$ we find that $R^1\pi_*\F$ and $R^1\pi_*\F(\Delta)$ are line bundles. As the map $R^1\pi_*\F \to  R^1\pi_*\F(\Delta)$ is a surjection it has a trivial kernel and we end up with the short exact sequence
\begin{equation}
  0\to  \O_{\tSop}\to E\to \pi_*(\F\otimes\omega_{\pi}^{\vee}|_{\Delta}) \to 0.
  \label{eq:exact_sequence_for_E}
\end{equation}
Introduce the notation
\begin{equation}
  \psi_1 := c_1(\pi_*(\omega_{\pi}|_{\Delta}))
  \label{eq:definition_of_psi1}
\end{equation}
It follows from~\eqref{eq:c_1_F} and~\eqref{eq:exact_sequence_for_E} that
\begin{equation}
  c_1(E) = \frac{\lambda}{4}- \frac{\psi_1}{2}.
  \label{eq:first_chern_class_of_E}
\end{equation}

Let us now consider the projective bundle
\begin{equation}
  P: = \mathbb PE\to \tSop.
  \label{eq:definition_of_P}
\end{equation}
We denote by $T\to P$ the tautological line bundle. Consider the fibered product $P\times_{\tSop} \tCop$ and let $P_1$ be the vanishing locus of the morphism $h: \pi_1^*T\to \pi_2^*\F(\Delta)$ (recall the definition of $E$ given in~\eqref{eq:DefOfE}), that is
\begin{equation}
  P_1 = \{(C, \eta, p, \sigma, q) \ \mid\ (C, \eta, p)\in \tSop,\quad \sigma\in \mathbb PH^0(C, \eta + p),\quad q\in\supp( \div \sigma)\}.
  \label{eq:DefOfP1}
\end{equation}
An easy calculation shows that $h$ has a zero of a first order along $P_1$, so that
\begin{equation}
  \div h = P_1.
  \label{eq:divisor_of_h_is_P1}
\end{equation}
The variety $P_1$ is a ``better version'' of $\tCop$ (cf.~\eqref{eq:DefOfP1} and~\eqref{eq:DefOftCop}). Indeed, the morphism $P_1\to \tCop$ is generically of degree $1$ because $|\eta+p|$ has no base points for a generic $(C,\eta,p)\in \tSop$, thus $P_1$ is birationally equivalent to $\tCop$. On the other hand the presence of ``$\sigma$'' in the description of points of $P_1$ affects well on the smoothness (cf.~proof of Lemma~\ref{lemma:tSopt_is_smooth}).

\subsection{Relation between $[\W]$ and $\lambda$.}

Recall that we defined $\W\subset \So$ as
\begin{equation*}
  \W = \text{close of}\left\{ (C, \eta)\in \So\ \mid\ \text{$|\eta + p|$ has a base point for some $p\leq \eta$} \right\}.
\end{equation*}
The set on the right-hand side is non-empty and the condition defining it is divisorial. It follows that $\W$ is either a non-empty divisor or coincides with $\So$. Let us show that the first case holds.

\begin{lemma}
  For any $g\geq 2$ there exists an odd spin curve $(C, \eta)$ and a point $p\in \supp(\eta)$ such that $h^0(C, \eta) = 1$ and $|\eta + p|$ does not have a base point. 
  \label{lemma:W_is_a_divisor}
\end{lemma}

\begin{proof}
  The proof goes by induction on $g$. If $g = 2$ then the statement is obvious, so let us suppose that $g\geq 3$ and there exists an odd spin curve $(D, \eta_D)$ of genus $g-1$ and $p\in D$ satisfying properties stated in the lemma. Let $y\in D\smm \supp(\eta_D)$ be a point that is not a ramification point of $|\eta_D + p|$. Pick a smooth elliptic curve $X$ with a marked point $x$ and an even theta characteristic $\eta_X$ on it and consider the odd spin curve given by
  \begin{equation*}
    (C, \eta) := (D\cup_y\mathbb P^1 \cup_x X,  \eta_D, \O(1), \eta_X).
  \end{equation*}
  Now let $(C_t, \eta_t)$ be a generic deformation of $(C, \eta)$ in $\cSo$ and $p_t\in \supp(\eta_t)$ be a family of points such that $p_0 = p\in D$. Suppose that for any $t$ there exists a $q_t\in \supp(\eta_t)$ that is a base point of $|\eta_t + p_t|$. Then $q_t$ cannot converge to $D$ since $|\eta_D + p|$ does not have base points, so $q_t$ converges to the unique point $q\in X$ such that $\eta_X = q-x$. But then a local analysis shows that the linear system $|q + x|$ has a base point which is clearly impossible. It follows that for a generic $t$ the linear system $|\eta_t + p_t|$ is base point free. Choosing $t$ small enough we can achieve small enough we can achieve the property $h^0(C_t, \eta_t) = 1$ also holds.
\end{proof}

We conclude that $\W$ does not coincide with all $\So$, so it is a divisor. In particular this implies that the projection from $P_1$ onto $\tCop$ is one-to-one over a generic $(C,\eta,p,q)\in \tCop$, since for any $p\in \supp(\eta)$ and $q\in C$ there exists a unique $\sigma\in |\eta + p|$ that passes through $q$. The locus where this projection is not one-to-one can be described by 
\begin{equation}
  \Wopt = \left\{ (C, \eta, p, q)\ \mid\ h^0(\eta + p - q) = 2 \right\}\subset \tCop
  \label{eq:def_of_Wopt}
\end{equation}
(recall that $h^0(\eta + p - q) = 2$ is equivalent to the fact that $q$ is a base point for $|\eta+p|$). The locus $\Wopt$ has codimension $2$, thus it follows that the map $P_1\to \tCop$ is the blow up of $\Wopt$. Let us use this observation to relate $[\W]$ and $\lambda$ in the rational Picard group of $\So$. Introduce the notation
\begin{equation}
  t := c_1(T).
  \label{eq:definition_of_t}
\end{equation}
It follows from~\eqref{eq:divisor_of_h_is_P1} and~\eqref{eq:c_1_F} that
\begin{equation}
  P_1 = \div h \equiv  \pi_2^*\Bigl(c_1(\F(\Delta))\Bigr) - \pi_1^*t \equiv \pi_2^*\left (\frac{1}{2}c_1(\omega_{\pi}) + \Delta\right ) + \frac{\lambda}{4}- \pi_1^*t
  \label{eq:P1_is_equal_to_difference}
\end{equation}
(recall that $\lambda$ always means the pullback of the Hodge class under the natural map to the moduli space $\cM$). Using the adjunction formula and the relation 
\begin{equation*}
  K_{P\times_{\tSop} \tCop} \equiv \pi_1^*(K_P - K_{\tSop}) + \pi_2^*K_{\tCop} \equiv \pi_1^*(2t-c_1(E)) + \pi_2^*K_{\tCop}
\end{equation*}
we deduce from~\eqref{eq:P1_is_equal_to_difference}:
\begin{equation}
  \begin{split}
    K_{P_1} & \equiv (K_{P\times_{\tSop} \tCop} + P_1)\cdot P_1 \\
    & \equiv \pi_2^* K_{\tCop}\cdot P_1 + \left( \pi_1^*(t-c_1(E)) + \pi_2^*\left (\frac{1}{2}c_1(\omega_{\pi}) + \Delta\right ) + \frac{\lambda}{4} \right)\cdot P_1
  \end{split}
  \label{eq:formula_for_K_P_1_via_adjunction}
\end{equation}
On the other hand we can compute $K_{P_1}$ using that the map $P_1\to \tCop$ is a blow up. Let $\eWopt\subset P_1$ be the exceptional divisor, that is, $\eWopt$ is the preimage of $\Wopt$, then
\begin{equation}
  K_{P_1} \equiv \pi_2^*K_{\tCop} \cdot P_1 + \eWopt.
  \label{eq:K_P_1_via_pi_2}
\end{equation}
It follows that
\begin{equation}
  \eWopt \equiv \left (\pi_1^*(t - c_1(E)) + \pi_2^*\left (\frac{1}{2} c_1(\omega_{\pi}) + \Delta\right) + \frac{\lambda}{4} \right )\cdot P_1.
  \label{eq:tildeW_via_K_P_1}
\end{equation}
The pushforward of $\eWopt$ under the projection $\tCop\to \So$ is equal to $\W$. To obtain a relation between $\W$ and $\lambda$ we intersect both sides of~\eqref{eq:tildeW_via_K_P_1} with a divisor in $P_1$ that is constructed below and push the result forward to $\So$. Let us first construct the divisor. We start with the divisor $\widehat \Delta\subset \tCop$ defined by the relation $\div (\pi^*\pi_* \F\to \F) = \Delta + \widehat \Delta$. Note that $\widehat \Delta$ can be described as
\begin{equation}
  \widehat \Delta =\{(C, \eta, p, q) \in \tCop\ \mid\ \eta\geq p+q\}.
  \label{eq:definition_of_widehat_Delta}
\end{equation}
The pullback of $\widehat \Delta$ to $P_1$ under the projection $P_1\to \tCop$ splits into the union $\widehat \Delta_1 \cup \eWopt$, where
\begin{equation}
  \widehat\Delta_1 = \left\{ (C, \eta, p, \sigma, q) \in P_1\ \mid\  \eta\geq p+q, \quad \sigma=\div\eta\subset |\eta + p| \right\}.
  \label{eq:Delta_1_as_set}
\end{equation}
Observe that
\begin{equation}
  (\pi_2)_*(\widehat \Delta_1\cdot \eWopt) = \Wopt,
  \label{eq:pushforward_of_eWopt_cdot_Delta_1_to_Cop}
\end{equation}
where the product is taken inside the Chow ring of $P_1$. Introduce a temporary notation for mappings:
\begin{equation}
  \xymatrix{
    P_1 \ar[d]^u  \ar@{^{(}->}[rr] && P\times_{\tSop} \tCop\ar[dll]^w \\
    \tSop \ar[d]^v\\
    \So
  }
\end{equation}
Using this notation and~\eqref{eq:tildeW_via_K_P_1} we write
\begin{equation}
  \begin{split}
    u_*(\eWopt\cdot \widehat \Delta_1) & \equiv u_* (\left( \left (\pi_1^*(t - c_1(E)) + \pi_2^*(\frac{1}{2} c_1(\omega_{\pi}) + \Delta)+ \frac{\lambda}{4}\right )\cdot P_1 \right) \cdot \widehat \Delta_1)\\
    & \equiv w_* ( \left (\pi_1^*(t - c_1(E)) + \pi_2^*(\frac{1}{2} c_1(\omega_{\pi}) + \Delta)+ \frac{\lambda}{4}\right )\cdot \widehat \Delta_1),
  \end{split}
  \label{eq:pushforward_of_eWopt_cdot_Delta_1_to_Sop}
\end{equation}
where the last product is taken inside the Chow ring of $P\times_{\tSop} \tCop$. In order to open the brackets we need the following

\begin{lemma}
  The following relations hold in $\Pic(\tSop)\otimes \mathbb Q$:
  \begin{align}
    & w_*((\pi_1^*c_1(E))\cdot \widehat\Delta_1) \equiv (g-2)\, c_1(E), \label{eq:c_1_E_times_Delta_1}\\
    & w_*\left( (\pi_1^*t)\cdot \widehat\Delta_1 \right) \equiv 0, \label{eq:t_times_Delta_1}\\
    & w_* \left ( (\pi_2^*\Delta)\cdot\widehat\Delta_1\right ) \equiv \Zop,  \label{eq:Delta_times_Delta_1} \\
    & w_*\left( (\pi_2^*c_1(\omega_{\pi})) \cdot \widehat \Delta_1 \right) \equiv v^*\mu - \psi_1,  \label{eq:omega_pi_times_Delta_1}
  \end{align}
  where $\Zop$ is defined by~\eqref{eq:definition_of_Zop} and $\mu$ is defined by~\eqref{eq:definition_of_mu}.
  \label{lemma:multiplication_by_Delta_1_over_smooth_curves}
\end{lemma}

We postpone the proof of this lemma till the end of the subsection, and now use it to open the brackets in~\eqref{eq:pushforward_of_eWopt_cdot_Delta_1_to_Cop}. From the description~\eqref{eq:Delta_1_as_set} one concludes that the map $\widehat\Delta_1\to\tSop$ has degree $g-2$. We obtain
\begin{equation*}
  u_*(\eWopt\cdot \widehat \Delta_1) \equiv -(g-2) c_1(E) + \frac{1}{2} \left( v^*\mu - \psi_1 \right) + \Zop + \frac{g-2}{4}\lambda.
\end{equation*}
Using~\eqref{eq:first_chern_class_of_E} we get
\begin{equation}
  u_*(\eWopt\cdot \widehat \Delta_1) \equiv \frac{g-3}{2}\,\psi_1 + \frac{1}{2} v^*\mu + \Zop.
  \label{eq:W_on_Sop_formula}
\end{equation}
Pushing this forward to $\So$ we obtain the following
\begin{prop}
  The following relation holds in $\Pic(\So)\otimes \mathbb Q$:
  \begin{equation*}
    \W \equiv \frac{g^2 + 11g - 6}{4}\, \lambda.
  \end{equation*}
  \label{prop:relation_between_W_and_lambda}
\end{prop}

\begin{proof}
  Note that $v_*u_*(\eWopt\cdot \widehat \Delta_1) \equiv 2\W$: it follows from the observation that the quantity $h^0(C,\eta+p-q)$ in~\eqref{eq:def_of_Wopt} is symmetric in $p$ and $q$, thus the map $\Wopt\to \W$ is of degree two. Applying $v_*$ to the right-hand side of~\eqref{eq:W_on_Sop_formula} and using that $v_*\psi_1 = \mu$ we get
  \begin{equation*}
    2\W\equiv (g-2)\mu + \Z.
  \end{equation*}
  By Lemma~\eqref{lemma:relation_between_mu_and_lambda} we have $\mu \equiv \frac{g+11}{2}\, \lambda$. Moreover, the formula from Theorem~\ref{thma:Farkas_formula} implies that $\Z \equiv (g+8)\,\lambda$ in the Chow ring of $\So$. Substituting this into the equality above we obtain the desired relation.
\end{proof}

\begin{proof}[Proof of Lemma~\ref{lemma:multiplication_by_Delta_1_over_smooth_curves}]
  The equality~\eqref{eq:c_1_E_times_Delta_1} follows from the fact that $w|_{\widehat\Delta_1}$ has degree $g-2$, while $\pi_1^*c_1(E)$ is a pullback of a divisor class from $\tSop$. To prove~\eqref{eq:t_times_Delta_1} we note that by~\eqref{eq:exact_sequence_for_E} the trivial line bundle $\O_{\tSop}$ injects into $E$, so that we obtain a section $\tSop\to P$, whose image we will denote by $D_{\mathrm{triv}}$. Notice that $T|_{D_{\mathrm{triv}}}$ is isomorphic to the trivial line bundle. The desctiption~\eqref{eq:Delta_1_as_set} of $\widehat\Delta_1$ implies that the image of $\widehat\Delta_1$ under the projection $\pi_1: P\times_{\tSop}\tSop\to P$ coincides with $D_{\mathrm{triv}}$ and the map $\pi_1\widehat\Delta_1\to D_{\mathrm{triv}}$ is $(g-2):1$. It follows that
  \begin{equation*}
    w_*\left( (\pi_1^* t)\cdot \widehat\Delta_1 \right) = (g-2)\,u_*\,c_1\left( T|_{D_F} \right) = 0.
  \end{equation*}

  Let $B\in A^1(\tCop)$ be arbitrary. Then we have $\pi_2^*B\cdot \eWopt = 0$, since $\eWopt$ is the exceptional divisor of $\pi_2|_{P_1}$. Using this observation we can write
  \begin{equation*}
    w_* \left ( (\pi_2^*B)\cdot\widehat\Delta_1\right ) = \pi_*(B\cdot \widehat\Delta),
  \end{equation*}
  where $\pi: \tCop\to \tSop$ is the forgetful map. Indeed, as $\widehat \Delta_1 \cup \eWopt$ is the pullback of $\widehat \Delta$ to $P_1$ and the restriction of $\pi_2$ to $P_1$ is of degree $1$ we can write
  \begin{equation*}
    w_* \left ( (\pi_2^*B)\cdot\widehat\Delta_1\right ) = w_*\left( (\pi_2^*B)\cdot\Bigl((\pi_2^*\widehat \Delta)\cdot P_1 - \eWopt\Bigr) \right) = u_*\left( (\pi_2|_{P_1})^*(B\cdot\widehat\Delta) \right) =  \pi_*(B\cdot \widehat\Delta).
  \end{equation*}

  To show~\eqref{eq:Delta_times_Delta_1} we substitute $B = \Delta$ and use the equality $\pi_*(\Delta\cdot \widehat\Delta) \equiv \Zop$ that can be checked easily.
  
  To show~\eqref{eq:omega_pi_times_Delta_1} we substitute $B = c_1(\omega_{\pi})$. We find that~\eqref{eq:omega_pi_times_Delta_1} follows from $\pi_*(c_1(\omega_{\pi})\cdot \widehat\Delta) \equiv v^*\mu - \psi_1$. Let $\phi: \tCop\to \tCo$ be the forgetful map, so that we have a commutative diagramm
  \begin{equation}
    \xymatrix{ \tCop \ar[rr]^{\phi} \ar[d]^{\pi} && \tCo \ar[d]^l \\ \tSop \ar[rr]^v && \tSo }
    \label{eq:commutative_diagram_Cop_Co_Sop_So}
  \end{equation}
  Note that $\omega_{\pi} \simeq \phi^*(\omega_l) $. We have $\Delta\cup\widehat \Delta = \phi^*(\tSop)$, so by definitions of $\mu$ (see~\eqref{eq:definition_of_mu}) and of $\psi_1$ (see~\eqref{eq:definition_of_psi1}) we find that
  \begin{equation*}
    \begin{split}
      \pi_*(c_1(\omega_{\pi})\cdot \widehat\Delta) & \equiv \pi_*( c_1(\omega_{\pi})\cdot (\phi^*(\tSop) - \Delta))\\
      & \equiv \pi_*\phi^*(c_1(\omega_l)\cdot\tSop) - \psi_1\\
      & \equiv v^*l_*(c_1(\omega_l)\cdot\tSop) - \psi_1\\
      & \equiv v^*\mu - \psi_1.
    \end{split}
  \end{equation*}

\end{proof}

\subsection{Relation between $\Cau$ and $\lambda$.}\label{subsec:relation_between_Cau_and_lambda}

The goal of this subsection is to prove the following

\begin{prop}
  The following relation holds in $\Pic(\So)\otimes \mathbb Q$:
  \begin{equation*}
    \Cau \equiv \frac{9g^2 + 179g -134}{2}\,\lambda.
  \end{equation*}
  \label{prop:relation_between_Cau_and_lambda}
\end{prop}

\noindent To prove Proposition~\ref{prop:relation_between_Cau_and_lambda} we will use the observation that the projection $P_1\to P$ is a branched cover: preimages of any element $(C, \eta, p, \sigma)\in P$ (where $\sigma\in |\eta + p|$) corresponds to different zeros of $\sigma$. The ramification divisor of this projection is therefore equal to the locus of those $\sigma$ that have multiple zeros. It has two components $B_1$ and $\tSopt$, where
\begin{align}
  & B_1 := \left\{ (C, \eta, p, \sigma, p)\in P_1\ \mid\ \div(\sigma)\geq 2p \right\}, \label{eq:definition_of_Delta_0}\\
  & \tSopt := \text{closure of } \left\{ (C, \eta, p, \sigma, q)\in P_1\ \mid\ \div(\sigma) = 2q + x_1 + \dots + x_{g-2},\quad p\neq q \right\}. \label{eq:definition_of_tSopt}
\end{align}
Note that
\begin{equation}
  B_1 \equiv \pi_2^*\Delta\cdot P_1.
  \label{eq:B_1_via_Delta}
\end{equation}
It is straightforward that $P_1\to P$ has a simple ramification along $B_1\cup \tSopt$. 

The projection $\varphi: \tSopt\to \tSop$ is a branched cover and its branched locus is given by two components $\BW$ and $\BCau$ lying above $\W$ and $\Cau$. Here
\begin{align}
  & \BW := \left\{ (C, \eta, p, \sigma, q) \in \tSopt\ \text{$q$ is a base point of $|\eta+p|$} \right\}, \label{eq:definition_of_BW}\\
  & \BCau := \left\{ (C, \eta, p, \sigma, q)\in \tSopt\ \mid\ \div(\sigma) = 3q + x_1 + \dots + x_{g-3} \right\}. \label{eq:definition_of_BCau}
\end{align}
As we will see (cf.~Corollary~\ref{cor:tSopt_to_tSop_simple_ramification}), the ramification is simple along these locuses. It follows that $K_{\tSopt} - \phi^* K_{\tSop} \equiv \BW + \BCau$ provided some smoothness of $\tSopt$. The pushforward of $\BCau$ to $\So$ is equal to $\Cau$, while the pushforward of $\BW$ is equal to $2\W$. It follows that if we will be able to relate the pushforward of $K_{\tSopt} - \phi^* K_{\tSop}$ and $\lambda$ in the Chow ring of $\So$ then we will obtain the desired relation between $\Cau$ and $\lambda$ (as we already know such a relation between $\W$ and $\lambda$ due to Proposition~\ref{prop:relation_between_W_and_lambda}). This is the plan and let us start with an appropriate smoothness statement for $\tSopt$:

\begin{lemma}
  The following properties of $\tSopt$ and $\BW$ and $\BCau$ holds:
  \begin{enumerate}
    \item The locus $\BCau$ is irreducible and we have $\div(\sigma) = 3q + x_1+\dots + x_{g-3}$ with $x_i\in C\smm\{q\}$ and $x_j\neq x_j$ if $i\neq j$ for a generic point $(C,\eta,p,\sigma,q)$ .\\
    \item The locus $\BCau$ is not a component of $\BW$.\\
    \item The variety $\tSopt$ is smooth.
  \end{enumerate}
  \label{lemma:tSopt_is_smooth}
\end{lemma}

As an immediate corollary we obtain the following statement:
\begin{cor}
  The ramification of the mapping $\tSopt\to \tSop$ is simple.
  \label{cor:tSopt_to_tSop_simple_ramification}
\end{cor}

We postpone the proof of this lemma to the end of the subsection. Recall the notation for the projection $w: P\times_{\tSop} \tCop\to \tSop$. Since $P\to \tSop$ is the projective bundle associated with the vector bundle $E$ we have $K_P\equiv K_{\tSop} + 2t - c_1(E)$. Using that $P_1\to P$ has a simple ramification along $B_1\cup\tSopt$,~\eqref{eq:B_1_via_Delta} and these observations we can write
\begin{equation}
  \tSopt\equiv K_{P_1} - (\pi_1^*K_P+ \pi_2^*\Delta)\cdot P_1 \equiv K_{P_1} - (\pi_1^*(2t - c_1(E)) + \pi_2^*\Delta + w^*K_{\tSop})\cdot P_1.
  \label{eq:tSopt_via_difference_of_kanonical_classes}
\end{equation}
Using the formula~\eqref{eq:formula_for_K_P_1_via_adjunction} and $K_{\tCop}\equiv\equiv  \pi^*K_{\tSop} + c_1(\omega_{\pi})$ we can write
\begin{equation}
  \begin{split}
    K_{P_1} & \equiv \left (\pi_2^* K_{\tCop} + \pi_1^*(t - c_1(E)) + \pi_2^*\left (\frac{1}{2} c_1(\omega_{\pi}) + \Delta\right) + \frac{\lambda}{4}\right )\cdot P_1\\
    &\equiv \left (w^* K_{\tSop} + \pi_1^*(t - c_1(E)) + \pi_2^*\left (\frac{3}{2} c_1(\omega_{\pi}) + \Delta\right) + \frac{\lambda}{4}\right )\cdot P_1
  \end{split}
  \label{eq:formula_for_K_P_1_via_adjunction_again}
\end{equation}
Substituting this expression for $K_{P_1}$ into~\eqref{eq:tSopt_via_difference_of_kanonical_classes} we find the following relation in the Chow ring of $P_1$:
\begin{equation}
  \tSopt \equiv \left( -\pi_1^*t + \frac{3}{2}\pi_2^*c_1(\omega_{\pi}) + \frac{\lambda}{4} \right)\cdot P_1.
  \label{eq:tSopt_via_other_classes}
\end{equation}
Lemma~\ref{lemma:tSopt_is_smooth} allows us to apply the adjunction formula to compute $K_{\tSopt}$. Combining~\eqref{eq:formula_for_K_P_1_via_adjunction_again} and~\eqref{eq:tSopt_via_other_classes} we find that
\begin{equation}
  \begin{split}
    K_{\tSopt} & \equiv (K_{P_1} + \tSopt)\cdot \tSopt\\
    & \equiv \left( w^*K_{\tSop} - \pi_1^*c_1(E) + \pi_2^*(3c_1(\omega_{\pi}) + \Delta) + \frac{\lambda}{2} \right)\cdot \tSopt
  \end{split}
  \label{eq:K_tSopt_via_adjunction}
\end{equation}
Recall that the forgetful mapping $\phi:\tSopt\to \tSop$ has a simple ramification along $\BW\cup \BCau$, so we have $\BCau \equiv K_{\tSopt} - w^* K_{\tSop}\cdot \tSopt - \BW$, where we replaced $\phi^*K_{\tSop}$ with $w^* K_{\tSop}\cdot \tSopt$. Thus the relation~\eqref{eq:K_tSopt_via_adjunction} implies 
\begin{equation}
  \begin{split}
    \phi_*\BCau & \equiv \phi_*(K_{\tSopt} - w^* K_{\tSop}\cdot \tSopt - \BW)\\
    & \equiv \phi_*\left( \left(-\pi_1^*c_1(E) + \pi_2^*(3\omega_{\pi} + \Delta) + \frac{\lambda}{2}\right)\cdot \tSopt - \BW \right) .
  \end{split}
  \label{eq:BCau_via_K_tSopt}
\end{equation}
Recall that $v_*\phi_*\BCau \equiv \Cau$ where $v: \tSop\to \So$. Thus it remains to relate the image of the right-hand side of~\eqref{eq:BCau_via_K_tSopt} under $v_*$ and $\lambda$. To do this let us introduce a divisor $\Zhat$ on $\tSop$ defined such that $v^*\Z \equiv \Zop + \Zhat$ where $\Zop$ was defined in~\eqref{eq:definition_of_Zop}. We can define $\Zhat$ alternatively as
\begin{equation}
  \Zhat := \left\{ (C, \eta, p) \in \tSop\ \mid\ \eta \geq p + 2q \text{ for some $q\in C$} \right\}.
  \label{eq:definition_of_Zhat}
\end{equation}

\begin{lemma}
  We have
  \begin{equation*}
    \phi_*( (\pi_2^*c_1(\omega_{\pi}))\cdot \tSopt) \equiv 4\Zop + 2\Zhat + 2\phi_*\BW - \frac{4g-3}{2}\lambda.
  \end{equation*}
  \label{lemma:t_cdot_tSopt}
\end{lemma}

\begin{proof}
  Recall that we have $c_1(\F) = \frac{1}{2}c_1(\omega_{\pi}) + \frac{\lambda}{4}$ in the Chow ring of $\tCop$ due to the relation~\eqref{eq:c_1_F}. On the other hand, we have $c_1(\F) \equiv \Delta + \widehat\Delta$ (cf.~the definition of $\tSop$ in Section~\ref{sec:spinors_with_marked_zero}) where $\widehat\Delta$ was introduced in~\eqref{eq:definition_of_widehat_Delta}. It follows that
  \begin{equation*}
    c_1(\omega_{\pi}) \equiv 2(\Delta + \widehat\Delta) - \frac{\lambda}{2}.
  \end{equation*}
  Note that $\phi_*(\pi^*\Delta\cdot \tSopt) = \Zop$ while $\phi_*(\pi^*\widehat\Delta\cdot \tSopt) = \phi_*( (\widehat\Delta_1 + \eWopt)\cdot \tSopt) = \Zop + \Zhat + \phi_*\BW$. As the projection $\tSopt\to \tSop$ has degree $4g-3$ (which is equal to the number of ramification points of $|\eta+p|$ not equal to $p$, cf.~\eqref{eq:definition_of_tSopt}), we have $\phi_*(\lambda\cdot \tSopt) = (4g-3)\lambda$. The lemma follows.
\end{proof}

Using this lemma and~\eqref{eq:first_chern_class_of_E} we can rewrite~\eqref{eq:BCau_via_K_tSopt} as
\begin{align}
  \phi_* \BCau & \equiv -(4g-3)c_1(E) + 13\Zop + 6\Zhat + 5\phi_*\BW - 2(4g-3)\,\lambda\\
  &\equiv \frac{4g-3}{2}\psi_1 + 13\Zop + 6\Zhat + 5\phi_*\BW - \frac{5(4g-3)}{4}\,\lambda
  \label{eq:Cau_on_tSop}
\end{align}
Now let us proof of Proposition~\ref{prop:relation_between_Cau_and_lambda}.

\begin{proof}[Proof of Proposition~\ref{prop:relation_between_Cau_and_lambda}]
  Note that $v_*\phi_*\BW \equiv 2\W$ in the rational Picard group of $\So$. Pushing~\eqref{eq:Cau_on_tSop} to $\So$ under $v: \tSop\to \So$ and using Proposition~\ref{prop:relation_between_W_and_lambda}, Lemma~\ref{lemma:relation_between_mu_and_lambda} and Theorem~\ref{thma:Farkas_formula} we get
  \begin{equation*}
    \begin{split}
      \Cau & \equiv \frac{4g-3}{2}\,\mu - \frac{5(4g-3)(g-1)}{4}\, \lambda + (6g-5)\Z + 10\W\\
      & \equiv \frac{9g^2 + 179 g - 134}{2} \,\lambda.
    \end{split}
  \end{equation*}
\end{proof}

\begin{proof}[Proof of Lemma~\ref{lemma:tSopt_is_smooth}]
  In this proof we again use the theory of translation surfaces as we did in the proof of Lemma~\ref{lm:tSop_is_smooth}. But in this cas we are dealing with meromorphic sections of spin bundles, so we need an analogous construction for meromorphic differentials which we describe now. Let $pr: \Mp\times_{\M}\Mp\to \Mp$ be the projection onto the first factor and let $p: \Mp\to \Mp\times_{\M}\Mp$ be the diagonal section. Consider the twisted Hodge vector bundle $\Ep:=pr_*\omega_{pr}(2p)$ Notice that fibers of $\Ep$ can be identified with spaces of Abelian differentials of the third kind with a double pole at the marked point $p$. Denote by $\mathbb P\Ep$ the corresponding projective bundle. Let $\mu = \{\mu_1,\dots, \mu_n\}$ be a set of integers such that $\mu_1\leq \dots\leq \mu_n$, $\mu_1\in \{-2\}\cup \mathbb Z_{>0}$, $\mu_2,\dots, \mu_n\in \mathbb Z_{>0}$ and $\sum_{i = 1}^n\mu_i = 2g-2$. Define $\Ad(\mu)\subset \mathbb P\Ep$ by
  \begin{equation*}
    \Ad(\mu) := \left\{(C,\omega)\in \mathbb P\Ep\ \mid\ \begin{array}{l} \div(\omega) = (2+\mu_1)p + \mu_2x_2 + \mu_3x_3+ \dots + \mu_nx_n,\\ \text{where $p$ is the marked point and $x_j\in C\smm\{p\}$, $x_i\neq x_j$} \end{array} \right\}
  \end{equation*}
  If one consider $\omega$ as a meromorphic differential then $(C,\omega)\in \Ad(\mu)$ iff $\mu$ is the set of multiplicities of $\omega$ at singular points.
  
  Introduce the partial order on partitions generated by relations $\mu'\prec \mu$ if $\mu'$ is obtained from $\mu$ by adding two of its elements. It is clear that
  \begin{equation}
    \text{closure of } \Ad(\mu) = \bigcup_{\mu'\preceq \nu} \Ad(\mu').
    \label{eq:closure_of_stratum}
  \end{equation}
  Now let us assume that all elements of $\mu$ are even and $(C,\omega)\in \Ad(\mu)$. Then $\frac{1}{2}\div(\omega) - p$ is a square root of $K_C$ in $\Pic(C)$, thus such an $\omega$ defines a spin structure. The parity of this spin structure is a topological invariant, so the locus $\Ad(\mu)$ has two components. Define by $\Ado(\mu)$ the component that corresponds to odd spin structures. The~\eqref{eq:closure_of_stratum} gives
  \begin{equation}
    \text{closure of } \Ado(\mu) = \bigcup_{\mu'\preceq \nu} \Ado(\mu').
    \label{eq:closure_of_odd_stratum}
  \end{equation}
  The following two facts are partial cases of main results of~\cite{ConComp} and~\cite{BOSSY}:
  \begin{enumerate}
    \item Let $\mu_+\subset \mu$ be the subset of positive elements of $\mu$. Then $\dim \Ado(\mu) = 2g + |\mu_+| - 2$.\\
    \item If $g\geq 4$ then the variety $\Ado(\mu)$ is connected for each $\mu$ satisfying $\min\{\mu_+\} = 2$. Moreover, $\Ado(\mu)$ is a smooth orbifold.
  \end{enumerate}
  
  Recall the definition of $P_1$ given in~\eqref{eq:DefOfP1}. Let $(C, \eta, p, \sigma, q)\in P_1$ and choose some isomorphism $\eta^{\otimes 2}\otimes \O_C(2p)\simeq K_C\otimes \O_C(2p)$. Let $\omega\in H^0(C,K_C+2p)$ be the image of $\sigma$ under this isomorphism. Then $(C,\omega)\in Cl\,\Ado(-2,2,2,\dots,2)$, where $Cl$ denotes the closure, so we get a map $P_1\to Cl\,\Ado(-2,2,2,\dots,2)$. This map has the degree $g$ and its image is dense in $\Ado(-2,2,2,\dots,2)$ (its complement is defined by conditions~\eqref{eq:no_automorphisms_and_h0=1}). The first advantage of this construction is that we can easily describe the image of $\BCau$ in terms of the partitions of zeros. Namely, the image of $\BCau$ is given by $Cl\,\Ado(-2,6,2,\dots,2)$ and the map $\BCau\to Cl\,\Ado(-2,6,2,\dots,2)$ is one-to-one over $\Ado(-2,6,2,\dots,2)$. Using the property (1) of $\Ado(\mu)$ and this observation we immediately get that $\BCau$ is irreducible.

  Now we will check that $\BCau$ is not a component of $\BW$. It is enough to give an example of $(C, \eta, p, \sigma, q)\in \BCau\smallsetminus \BW$. Consider a generic odd spin curve $(X, \eta_X)$ of genus $g-1$, pick a point $p\in \supp(\eta_X)$ and a ramification point $x\notin\supp(\eta_X)$ of $|\eta_X + p|$. Then we take a smooth elliptic curve $Y$ with an even theta characteristic $\eta_Y$ and $y,q\in Y$ such that $3q-3y\equiv \eta_Y$ but $q-y\not\equiv \eta_Y$. Then it is easy to see that the marked spin curve $(X\cup_x \mathbb P^1\cup_y Y, \eta_X, \O(1), \eta_Y, p,q)$ can be obtained as a limit of elements from $\BCau$ but cannot be obtained as a limit of elements of $\BW$ (otherwise the equation $q-y\equiv \eta_Y$ must hold).

  Now let us show that $\tSopt$ is smooth. Note that the image of $\tSopt$ in $Cl\,\Ado(-2,2,2,\dots,2)$ is given by $Cl\,\Ado(-2,4,2,\dots,2)$. Given a point $(C, \eta, p, \sigma, q)\in P_1$ we denote by $\omega_{\sigma}\in H^0(C,K_C+2p)$ any section that corresponds to $\sigma$. We will consider $\omega_{\sigma}$ as meromorphic differential. Assume that $\omega_{\sigma} = f(z)^2\,dz$ in some local coordinate $z$ near $q$ such that $z(q)=0$. Then it is easy to see that the map $(C, \eta, p, \sigma, q)\mapsto f'(0)$ can be extended in a small neighborhood of $(C,\eta,p,\sigma,q)$ in $\tSopt$ to a holomorphic function whose zero set coinside with $\tSopt$ in this neighborhood. Denote this function by $\Phi(\cdot)$. Our goal is to show that $\nabla\Phi$ does not vanish at $(C,\eta,p,\sigma,q)\in \tSopt$. Let us choose the coordinate $z$ such that $\omega_{\sigma} = z^{2k}\,dz$. Assume for simplicity that $p\neq q$ (the case $p=q$ can be made similarly). For a small $\eps>0$ and let $U_{\eps}\subset C$ be the closed neighborhood of $q$ defined by the condition $|z| \leq \eps$ and let $C_{\eps} = C\smallsetminus U_{\eps}$. For each $t\in \mathbb C$ consider the differential
  \begin{equation*}
    \nu_t := (x^k - tx - t^k + t^2)^2\,dx
  \end{equation*}
  on the complex plane with the coordinate $x$. Define the function $F$ on $\{(x,t)\ \mid\ \eps < |x| < 3\eps,\ |t|\ll \eps\}$ by the condition $F(x,0) = x$ and $d/dx\,(F(x,t)^{2k+1}) = (2k+1) (x^k - tx - t^k + t^2)^2$. For each small $t$ glue $C_{2\eps}$ with $\{x\in \mathbb C\ \mid\ |x| < 3\eps\}$ such that the condition $z = F(x,t)$ is satisfied. Denote the obtained Riemann surface by $C_t$. Differentials $\omega_{\sigma}$ and $\nu_t$ coincide on the glued part so that we obtain a global differential $\omega_t$ on $C_t$. Zeros and poles of $\omega_t$ corresponds to zeros and poles of $\omega_{\sigma}$ on $C_{2\eps}$ and to zeros of $\nu_t$ on $\{x\in \mathbb C\ \mid\ |x| < 3\eps\}$. Therefore all singularities of $\omega_t$ are even and $\eta_t:=\frac{1}{2}\div\omega_t$ is a spin structure. Since $\eta_0 = \eta$ is odd we see that $\eta_t$ is odd for all $t$ where defined. Let $\sigma_t\in \mathbb PH^0(C_t, \eta_t + p)$ be the section defined by the square root of $\omega_t$. The section $\sigma_t$ vanishes at $x=t$, so that we obtain a one-dimenstional family in $P_1$ given by the mapping $t\mapsto (C_t,\eta_t,p, \sigma_t, q_t)$ where $q_t\in C_t$ is the point corresponding to $x = t$. Then we have $\Phi(C_t,\eta_t,p,\sigma_t,q_t) = kt^{k-1} - t$. It follows that the gradient of $\Phi$ does not vanish at $t = 0$ and we are done.

\end{proof}

\begin{cor}
  The subvariety $\pi_2(\tSopt)\subset \tCop$ is smooth.
  \label{cor:cor_tSopt_is_smooth}
\end{cor}

\begin{proof}
  The projection $\pi_2:P_1\to \tCop$ is one-to-one outside $\eWopt$, thus $\pi_2(\tSopt)$ is smooth outside $\pi(\BW) = \Wopt$. We need to show that $\pi_2(\tSopt)$ is smooth at any point of $\Wopt$. To do this pick a $(C, \eta, p, q)\in \Wopt$ and let $\iota: \mathbb PH^0(C, \eta + p)\hookrightarrow P_1$ be the fiber over $(C, \eta, p,q)$. It is enough to show that $\iota^*\O_{P_1}(\tSopt) \simeq \O(1)$. Recall that $\tSopt + B_1$ is the ramification divisor of the map $P_1\to P$ and $B_1 \equiv P_1\cdot\pi_2^*\Delta$. As $(C, \eta, p, q)\in \Wopt$ the point $q$ is the base point of $|\eta+p|$ and thus $q\neq p$. It follows that $\iota^*\O_{P_1}(B_1) = 0$, since $\pi_2^{-1}(\Delta)\cap \iota(\mathbb PH^0(C, \eta + p)) = \varnothing$, and
  \begin{equation*}
    \iota^*\O_{P_1}(\tSopt) = \iota^*\O_{P_1}(\tSopt + B_1) = \iota^*\O_{P_1}(K_{P_1} - \pi_1*K_P\cdot P_1).
  \end{equation*}
  Recall~\eqref{eq:P1_is_equal_to_difference}:
  \begin{equation*}
    P_1 \equiv \pi_2^*c_1(\F(\Delta)) - \pi_1^*t.
  \end{equation*}
  Using this equality and the adjunction formula we get
  \begin{equation*}
    K_{P_1} - \pi_1^*K_P\cdot P_1 \equiv (K_{P\times_{\tSop}\tCop} + P_1 - \pi_1^*K_P)\cdot P_1 \equiv (-\pi_1^*t + \pi_2^*(\omega_{\pi} + c_1(\F(\Delta))))\cdot P_1
  \end{equation*}
  It follows that
  \begin{equation*}
    \iota^*\O_{P_1}(\tSopt) \equiv \iota^*\O_{P_1}\Bigl((-\pi_1^*t + \pi_2^*(\omega_{\pi} + c_1(\F(\Delta))))\cdot P_1\Bigr) \equiv \iota^*\O_{P_1}(-P_1\cdot \pi_1^*t) \equiv \O(1).
  \end{equation*}
\end{proof}

We finish this section with the following computation that we use in Section~\ref{sec:test_curves}.

\begin{lemma}
  Introduce the locus
  \begin{equation*}
    \X := \{(C, \eta, p, q)\in \tCop\ \mid\ \eta \geq p + 2q\}.
  \end{equation*}
  Then $\X$ is irreducible and we have
  \begin{equation*}
    \pi_2(\tSopt)\cdot \widehat\Delta = \X + \Wopt
  \end{equation*}
  where $\Wopt$ was defined in~\eqref{eq:def_of_Wopt} and $\widehat \Delta$ in~\eqref{eq:definition_of_widehat_Delta}.
  \label{lemma:Zg_cdot}
\end{lemma}

\begin{proof}
 It is straightforward that
  \begin{equation*}
    \pi_2(\tSopt)\cap\widehat\Delta = \X \cup \Wopt.
  \end{equation*}
  Thus, it remains to show that the intersection above is generically transversal. Let $(C, \eta, p, q)\in \X \cup \Wopt$. Assume first that $(C, \eta, p, q)\notin \X$. Then the projection $\pi_2(\tSopt)\to \tSop$ has a ramification at $(C,\eta, p, q)$ while $\widehat\Delta\to \tSop$ is unbranched at it, which implies the transversality of the intersection. Assume now that $(C,\eta, p,q)\in \X\smallsetminus \Wopt$. It follows that $\widehat\Delta\to \tSop$ has a branhing at this point while $\pi_2(\tSopt)$ has no ramification, therefore the intersection is again transversal. It remains to show that $\dim(\X\cap \Wopt) < \dim \X = \dim \Wopt$. Assume that $\X' = \X\cap \Wopt$ has the same dimension as $\X$. Since $\Z$ is irreducible (cf. the proof of Lemma~\ref{lm:tSop_is_smooth}) and the forgetful map $\X'\to\Z$ is of finite degree the image of $\X'$ coincides with $\Z\cap \tSo$. Let $(X,\eta_X)$ be a generic odd spin curve of genus $g-2$ and $x\in X$ be a point that is not a ramification point of $|\eta_X+p|$ for any $p\in \supp(\eta_X)$. Let $(Y,\eta_Y)$ be an even spin curve of genus $2$ such that $T_{\eta_Y}$ is smooth. Let $(y,q)\in T_{\eta_Y}$ be such that $\eta_Y + y\geq 2q$. Then it is straightforward (see~\cite[Section~5]{FARo}) that the odd spin curve $(X\cup_x E\cup_y Y, \eta_X,\O_E(1), \eta_Y)$ is a limit of a family of spin curves $(C_t,\eta_t)$ from $\Z$ and double points $q_t$ of $\eta_t$ converge to $q\in Y$. But as the image of $\X'$ coincides with $\Z\cap \tSo$ we can find a family of points $p_t\in \supp(\eta_t)\smallsetminus\{q_t\}$ such that $|\eta_t + p_t|$ has a base point at $q_t$. Note that points from $p_t$ must converge to a point $p\in \supp(\eta_X)$. Let $\sigma_X\in \mathbb PH^0(X, \eta_X\otimes \O_X(p + 2x)), \sigma_X\in \mathbb PH^0(Y, \eta_Y\otimes \O_Y( (g-2)y ))$ be aspects of $|\eta_t+p_t|$ (see~\cite{LimLin} for details). We should have $\ord_x\sigma_X + \ord_y\sigma_Y = g-1$. As $X$ is generic the maximal order of $\sigma_X$ at $x$ is $g-1$. It follows that $q$ should be a base point of $|\eta_Y + 2y|$, or, equivalently, $h^0(Y,\eta_Y + 2y-q) = 2$. It follows by Riemann-Roch that $h^0(Y,\eta_Y + q - 2y) = 1$ i.e. $\eta_Y + q\geq 2y$. But then the inequality $\eta_Y+y\geq 2q$ cannot be satisfied because $T_{\eta_Y}$ is smooth be the choice of $(Y,\eta_Y)$. Thus we obtain a contradiction.
\end{proof}

\section{Intersection with test curves.}\label{sec:test_curves}

The goal of this section is to calculate the boundary contribution in the formulas for $[\cCau]$ and $[\cW]$ and thus the proof Theorem~\ref{thm:caustic_and_base_point_divisors_formula}. Let us write
\begin{equation}
  \begin{split}
    & \cCau \equiv l^c \cdot \lambda - a_0^c \cdot \alpha_0 - b_0^c - \sum_{j = 1}^{g-1}a_j^c \cdot \alpha_j,\\
    & \cW \equiv l^w \cdot \lambda - a_0^w \cdot \alpha_0 - b_0^w - \sum_{j = 1}^{g-1}a_j^w \cdot \alpha_j\\
  \end{split}
  \label{eq:notation_for_coefficietns}
\end{equation}
in $\Pic(\cSo)\otimes \mathbb Q$. By Proposition~\ref{prop:relation_between_W_and_lambda} and Proposition~\ref{prop:relation_between_Cau_and_lambda} we have
\begin{equation}
  \begin{split}
    & l^c =  \frac{9g^2 + 179g -134}{2},\\
    & l^w = \frac{g^2 + 11g - 6}{4}.
  \end{split}
  \label{eq:coefficients_before_lambda}
\end{equation}
The other coefficients will be computed using intersections with test families. First, we compute the coefficients at $\alpha_j$ for $j>0$ and then finish the section computing coefficients at $\alpha_0$ and $\beta_0$.

\subsection{Universal deformation of a spin curve.}\label{subsec:universal_family}

Let us briefly recall how a universal deformation of a spin curve was constructed in~\cite{COR}. A smooth family of nodal curves $\pi: X\to B$ together with a line bundle $\L\to X$ and a homomorphism $f: \L^{\otimes 2}\to \omega_{\pi}$ is called a family of spin curves if each fiber $C$ of $\pi$ is a quasi-stable curve (see Section~\ref{sec:moduli_of_spin_curves}) and the triple $(C, \L|_C, f)$ is a spin curve. A morphism of two families of spin curves is defined in an obvious way. Given a spin curve $(C, \eta, \beta)$ a universal deformation of it is constructed as follows. Let $C'$ be the curve obtained from $C$ by contracting all exceptional components. Denote by $p_1,\dots, p_n\in C'$ the nodes that correspond to exceptional components in $C$. Let $B$ be a polydisc in $\mathbb C^{3g-3}$ such that $B/\Aut(C')$ is isomorphic to a small neighborhood of $C'$ in $\cM$ and let $\pi': X'\to B$ be the pullback of the universal curve. Assume moreover that coordinates $z_1,\dots,z_{3g-3}$ in $B$ are chosen such that $\{z_j = 0\}$ is the locus where the node $p_j$ persists. Let $B$ be an another polydisc (say, in coordinated $w_1,\dots, w_{3g-3}$) and $\varphi: B\to B'$ be the branched cover given by equations $z_j = w_j^2, j = 1,\dots, n$, $z_j = w_j, j = n+1,\dots, 3g-3$. For each $j = 1,\dots, n$ the variety $\varphi^* X'$ has a codimension $2$ subvariety consisting entirely of nodes and corresponding to $p_j$. Blowing up these subvarieties yields a family $\pi: X\to B$ with $X$ smooth and $C$ as a central fiber. Denote by $\CMcal E_j\subset X$ the blow up divisor corresponding to $p_j$. Set $\xi = \omega_{\pi}(-\sum \CMcal E_j)$. Then it can be shown (see~\cite{COR}) that, altering the idetification of $C$ with the central fiber by an isomorphism that is identical outside exceptional components, we can achieve that $\eta^{\otimes 2}$ is isomorphic to the restriction of $\xi$ to the central fiber. Shrinking $B$ if necessary we may then extend $\eta$ to a bundle $\L\to X$ and the isomorphism between $\eta^{\otimes 2} $ and $\xi|_C$ to the isomorphism between $\L^{\otimes 2}$ and $\xi$. Composing this with the natural monomorphism $\xi\to\omega_{\pi}$ we finally get a family of spin curves $(\pi: X\to B, \L, f)$ such that $(C, \eta, \beta)$ is isomorphic to the central fiber of this family. This family is called {\it universal deformation} of $(C, \eta, \beta)$ due to the following proposition (see~\cite[Proposition (4.6)]{COR}):

\begin{prop}[\cite{COR}]
  Let $(\pi_1: X_1\to T_1, \L_1, f_1)$ be a family of spin curves and $t\in T_1$ be a point. Let $(X\to B, \L, f)$ be the universal deformation of $(\pi_1^{-1}(t), \L_{\pi_1^{-1}(t)}, f_1)$. Then for any small neighborhood $B_1$ of $t$ there exists a unique morphism $(\pi_1^{-1}(B_1)\to B_1, \L_1, f_1)\to (X\to B, \pi, f)$ of families that identifies the central fiber of $(X\to B, \L, f)$ with $(\pi_1^{-1}(t), \L_{\pi_1^{-1}(t)}, f_1)$.
  \label{prop:universal_deformation_of_a_spin_curve_is_universal}
\end{prop}

\begin{rem}
  Let us denote by $\scSo$ the moduli {\it stack} of odd spin curves. It follows from Proposition~\ref{prop:universal_deformation_of_a_spin_curve_is_universal} that the morphism $\scSo\to \cSo$ is simply ramified over the divisor in $\cSo$ consisting of spin curves containing exceptional components.
  \label{rem:stack_to_space_is_ramified}
\end{rem}

\subsection{Test families for $\alpha_j$ for $j>0$.}\label{subsec:test_families_definition}

To compute the coefficient at $\alpha_j$ for $1\leq j\leq g-2$ we will use the standard covering families $G_j\subset A_j$ constructed as follows. Fix a generic odd spin curve $(X,\eta_X)$ of genus $j$ together with a generic even spin curve $(Y, \eta_Y)$ of genus $g-j$. Pick generic points $x\in X$ and $y\in Y$. Then we form an odd spin curve of genus $g$ by considering $(X\cup_x E\cup_y Y, \eta = \eta_C, \O_E(1), \eta_D)$ where $E$ is a rational component. Varying $y$ we obtain the family $G_j\subset A_j$. So,
\begin{equation}
  G_j := \left\{ X\cup_x E\cup_y Y, \eta = \eta_C, \O_E(1), \eta_D \right\}_{y\in Y}.
  \label{eq:definition_of_Gj}
\end{equation}
Classically, we have 
\begin{equation}
  \begin{split}
    & G_j\cdot \lambda = G_j\cdot \beta_0 = G_j\cdot \alpha_k = 0,\qquad k\neq j,\\
    & G_j\cdot \alpha_j = -2(g-j-1).
  \end{split} 
  \label{eq:Gj_cdot_basis}
\end{equation}

\begin{figure}[h]
  \begin{center}
    \includegraphics[width = 0.45\textwidth]{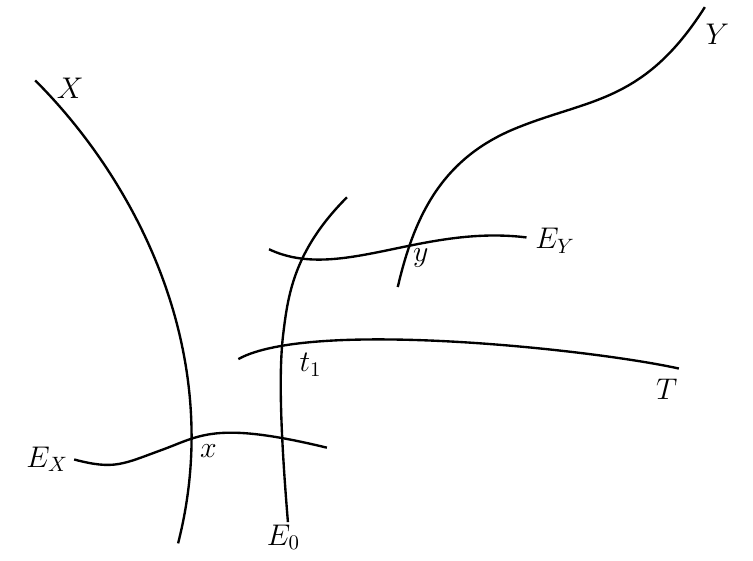}
  \end{center}
  \caption{Special fiber of $H_j$}
  \label{fig:special_fiber}
\end{figure}

To compute the coefficient at $\alpha_{g-1}$ we will use a more delicate construction of a test family. Let $2\leq j\leq g-1$ and consider a generic even spin curve $(X, \eta_X)$ of genus $j-1$ and a generic even spin curve curve $(Y, \eta_Y)$ of genus $g-j$. Let $T$ be a smooth elliptic curve with an odd theta characteristic $\eta_T$ and let $x\in X, y\in Y$ and $t_1,t_2\in T$ be some points; we assume that $x\in X$ and $y\in Y$ are chosen generically. For each pair $t_1\neq t_2$ we consider an odd spin curve
\begin{equation}
  (C_{t_1,t_2}, \eta) = (X\cup_x E_X\cup_{t_1} T\cup_{t_2} E_Y\cup_y Y, \eta_X, \O_{E_X}(1), \eta_T, \O_{E_Y}(1), \eta_Y)
  \label{eq:general_fiber_of_H}
\end{equation}
where $E_X, E_Y$ are rational components. When $t_1$ and $t_2$ collide we obtain a spin curve $(C_{t_1 = t_2}, \eta)$ that is schematically drawn on Figure~\ref{fig:special_fiber}, where $E_0$ is a rational component.

Varying the pair $(t_1,t_2)$ we obtain a family of odd spin curve that we denote by $H_j$. Note that $H_j$ is one-dimensional since the point in $H_j$ depends only on the difference $t_1-t_2$ of two points on the elliptic curve $T$, though we will assume that the point $t_1$ is fixed and the point $t_2$ moves. A straightforward computation shows that
\begin{equation}
  \begin{split}
    & H_j\cdot \lambda = H_j \cdot \beta_0 = 0,\qquad H_j \cdot\alpha_i = 0,\quad i \neq 1,g-j+1,j\\
    & H_j\cdot \alpha_j = H_j\cdot \alpha_{g-j+1} = -1,\\
    & H_j\cdot \alpha_1 = 1.
  \end{split}
  \label{eq:intersection_with_H}
\end{equation}

Notice that $H_j$ is homotopic to $H_{g-j+1}$, therefore we get $[g/2]-1$ equations on coefficients computing intersection numbers with $H_j$. We need only one of them to compute the coefficients at $\alpha_{g-1}$, thus the other $[g/2]-2$ can be used as a verification of our formulas.

\subsection{Intersection of $G_1,\dots, G_{g-2}$ with caustic and fixed point divisors.}\label{subsec:G_cap_Cau_and_W}

Fix some $j:\ 1\leq j\leq g-2$ till the end of the subsection. Before we start computing $G_j \cdot \cW$ and $G_j\cdot \cCau$ let us note that it is not hard to describe sets $G_j\cap \cW$ and $G_j\cap \cCau$ studying corresponding limits linear series (see Remark~\ref{rem:interpretaion_via_limit_series}). But to compute intersection numbers $G_j\cdot \cW$ and $G_j\cdot \cCau$ one needs to prove transversality (or determine the multiplicity) of the intersection and then to compute the size of the intersection. To handle this we will study constructions made in Section~\ref{sec:construction_over_the_smooth_part} locally near $G_j$. This allows us to reduce our questions to intersection theory of one-dimensional correspondences. As a result we reduce the question of computing $G_j\cdot \cW$ and $G_j\cdot \cCau$ to some computation on a product surface of type $Y\times\text{[some curve]}$.

The family $G_j$ is an image of the map $Y\to \cSo$ that sends $y\in Y$ to the moduli of the curve $(X\cup_x E\cup_y Y, \eta_X, \eta_Y)$. Let $\mathbb D\subset \mathbb C^{3g-4}$ be a small polydisc and set $B = Y\times \mathbb D$. Then we can extend $Y\to \cSo$ to a map $\CMcal{M}: B\to \cSo $ such that for each $y\in Y$ there exists a neighborhood $U\subset Y$ of $y$ such that the pullback of the universal spin curve under to $U\times \mathbb D$ is a universal deformation of $(X\cup_x E\cup_y Y, \eta_X, \eta_Y)$. Denote by $l : \C\to B$ the pullback of the universal spin curve under $\CMcal M$.

Let $\C\to \cCo$ be the natural map and let by $B_1\subset \C$ denote the closure of the preimage of $\tSop$. Let us describe connected components of $B_1$. Let $\eta_X \simeq \O_X(p_{X,1} + \dots + p_{X,j-1})$. We may assume that all $p_{X,i}$ are distinct as $X$ was chosen generic. Then $l^{-1}(Y\times \{0\})\cap B_1 = Y_1\sqcup \dots\sqcup Y_j$ where
\begin{align*}
  &Y_i = \{ (X\cup_x E\cup_y Y, \eta_X , \eta_Y, p_{X,i})\in \C \}_{y\in Y}, \qquad i = 1,\dots, j-1,\\
  &Y_j = \{ (X\cup_x E\cup_y Y, \eta_X , \eta_Y, p)\in \C \}_{(y,p)\in T_{\eta_Y}}.
\end{align*}
where $T_{\eta_Y} = \{(y,p)\in Y^2\ \mid\ \eta_Y + y \geq p\}$ is the Scorza curve associated with $\eta_Y$. Denote by $B_{X,i}$ the connected component of $B_1$ that contains $Y_i$ for each $i = 1,\dots, j-1$ and by $B_Y$ the connected component of $B_1$ that contains $Y_j$. Shrinking $\mathbb D$ if necessary we can achieve $B_1 = B_{X,1}\sqcup\dots\sqcup B_{X,j-1}\sqcup B_Y$.

Let $\C_1 = \C\times_B B_1$. Consider the following codimension $2$ locuses of $\C_1$:
\begin{align}
  & D_{\Cau} := \text{closure of } \{(C,\eta, p, q)\in \C_1\ \mid\ C\text{ is smooth and } \eta + p \geq 3q,\ q\neq p\},\\
  & D_{\W} := \text{closure of } \{(C,\eta, p,q)\in \C_1\ \mid\ C\text{ is smooth and } h^0(C, \eta + p - q) = 2\},\\
  & D_{\Z} := \text{closure of } \{(C,\eta, p,q)\in \C_1\ \mid\ C\text{ is smooth and } \eta \geq p+2q\}.
  \label{eq:def_of_DSubscript}
\end{align}
Let $\C_{X,i}\to B_{X, i}$ and $\C_Y\to B_Y$ be the restriction of $\C_1\to B_1$ to $B_{X,i}$ and $B_Y$. Introduce the following notation:
\begin{align*}
  & \iota_{X,i, X}: Y\times X\hookrightarrow \C_{X,i}\qquad \iota_{X,i,X}(y,q) = (X\cup_x E\cup_y Y, p_{X,i}, q),\\
  & \iota_{X,i, Y}: Y\times Y\hookrightarrow \C_{X,i}\qquad \iota_{X,i,X}(y,q) = (X\cup_x E\cup_y Y, p_{X,i}, q),\\
  & \iota_{Y, X}: T_{\eta_Y}\times X\hookrightarrow \C_{X,i}\qquad \iota_{Y,X}( (y,p),q) = (X\cup_x E\cup_y Y, p, q),\\
  & \iota_{Y, Y}: T_{\eta_Y}\times Y\hookrightarrow \C_{X,i}\qquad \iota_{Y,Y}( (y,p),q) = (X\cup_x E\cup_y Y, p, q).
\end{align*}

\begin{lemma}
  We have
  \begin{align}
    & \cCau \cdot G_j = \sum_{i = 1}^{j-1} (\deg\iota_{X,i,X}^*D_{\Cau} + \deg\iota_{X,i,Y}^*D_{\Cau}) + \deg\iota_{Y,X}^*D_{\Cau} + \deg\iota_{Y,Y}^*D_{\Cau},  \label{eq:DCau_instead_of_Cau_D}\\
    & \cW \cdot G_j = \frac{1}{2}\left( \sum_{i = 1}^{j-1} (\deg\iota_{X,i,X}^*D_{\W} + \deg\iota_{X,i,Y}^*D_{\W}) + \deg\iota_{Y,X}^*D_{\W} + \deg\iota_{Y,Y}^*D_{\W} \right),  \label{eq:DW_instead_of_W}\\
    & (g-3)\cZ \cdot G_j = \sum_{i = 1}^{j-1} (\deg\iota_{X,i,X}^*D_{\Z} + \deg\iota_{X,i,Y}^*D_{\Z}) + \deg\iota_{Y,X}^*D_{\Z} + \deg\iota_{Y,Y}^*D_{\Z}.  \label{eq:DZ_instead_of_Z}
  \end{align}

  \label{lemma:Gj_Ddivisor_instead_divisor}
\end{lemma}

Lemma~\ref{lemma:Gj_Ddivisor_instead_divisor} allows to reduce the problem of computing intersection numbers with $G_j$ to the computation of intersections between the varieties defined in~\eqref{eq:def_of_DSubscript} and two-dimensional families embedded into the boundary of $\C_1$. Notice that $D_{\Cau}$, $D_{\W}$ and $D_{\Z}$ corresponds to the varieties $\BCau$, $\BW$ and $\X$ introduced in Section~\ref{sec:construction_over_the_smooth_part}. Representing classes of $D_{\Cau}$, $D_{\W}$ and $D_{\Z}$ via intersections of divisors as we did for $\BCau$, $\BW$ and $\X$ in Section~\ref{sec:construction_over_the_smooth_part} we will compute the right-hand side of~\eqref{eq:DCau_instead_of_Cau_D}~-~\eqref{eq:DZ_instead_of_Z}.

Before we prove Lemma~\ref{lemma:Gj_Ddivisor_instead_divisor} let us make some observations. Denote by $\mathrm{pr}: \C_1\to B$ the forgetful map. Notice that $\mathrm{pr}_*D_{\Cau} \equiv \CMcal M^*\Cau\in A^1(B)$ , $\mathrm{pr}_* D_{\W} \equiv 2\CMcal M^*\W\in A^1(B)$ and $\mathrm{pr}_* D_{\Z} \equiv (g-3)\CMcal M^*\Z\in A^1(B)$. As $\CMcal{M}(Y\times \{0\}) = G_j$, it follows that $\cCau\cdot G_j = D_{\Cau}\cdot \mathrm{pr}^*(Y\times\{0\})$ and the same for the other two divisor classes. As the sum on the right-hand side in Lemma~\ref{lemma:Gj_Ddivisor_instead_divisor} is the intersection with a part of $\mathrm{pr}^*(Y\times \{0\})$ it follows that we need to show that $D_{\Cau}$, $D_{\W}$ and $D_{\Z}$ do not intersect the other part. Namely, we need to show that $D_{\Cau}$, $D_{\W}$ and $D_{\Z}$ do not intersect the exceptional component of any fiber of $\C_1$ over $Y_i\subset B_1$ for any $i = 1,\dots,j$. The fact that $D_{\W}$ and $D_{\Z}$ do not intersect exceptional components follows from the fact that zeros of sections of odd spin bundles cannot specialize to the exceptional component if a curve approach $X\cup_x E\cup_y Y,\eta_X,\O_E(1), \eta_Y$ if $x$ is chosen generic. Thus we are left to show the claim for $D_{\Cau}$. We will prove this by showing that the intersection number of $D_{\Cau}$ with this family of exceptional components is zero. In order to do this let us introduce several more varieties and represent $D_{\Cau}$ as an intersection of divisors in $\C_1$. Set
\begin{align}
  & \S_{\Sopt} := \text{closure of } \{(C,\eta, p, q)\in \C_1\ \mid\ C\text{ is smooth and } \eta + p \geq 2q,\ q\neq p\},\\
  & \S_{\dhat} := \text{closure of } \{(C,\eta, p, q)\in \C_1\ \mid\ C\text{ is smooth and } \eta \geq p+q\}.
  \label{eq:def_of_SSubscript}
\end{align}
Notice that the variety $\S_{\Sopt}$ is the pullback of the variety $\pi_2(\tSopt)$ (see Section~\ref{subsec:relation_between_Cau_and_lambda} and in particular Corollary~\ref{cor:cor_tSopt_is_smooth}) while $\S_{\dhat}$ is the pullback of $\widehat\Delta$ (see~\eqref{eq:definition_of_widehat_Delta}). The forgetful map $\S_{\Sopt}\to B_1$ fails to be a branched cover over a locus in the boundary of $B_1$ which we describe now:

\begin{lemma}
  Consider the locus $D\subset \C_1$ that parametrizes the following set of marked spin curves:
  \begin{equation*}
    \{ (X_1\cup_{x_1} E\cup_{y_1} Y_1, \eta^-_{X_1}, \eta_{Y_1}^+, p, q)\in \C_1\ \mid\ p\in Y_1,\ \eta^+_{Y_1} + p \geq 2y_1,\ q\in X_1\cup E \}.
  \end{equation*}
  Then $D$ is a divisor in $\S_{\Sopt}$ and the forgetful map $\S_{\Sopt}\to B_1$ is a branched cover outside $D$. As a divisor $D$ splits into sum of two components $D_1+D_2$ where $D_1$ corresponds to the case when $q\in X_1$ and $D_2$ corresponds to the case when $q\in E$.
  \label{lemma:Dfib_is_fiber}
\end{lemma}

\begin{proof}
  Assume that $(X_1\cup_{x_1} E\cup_{y_1} Y_1, \eta^-_{X_1}, \eta_{Y_1}^+, p, q)\in \S_{\Sopt}$ is a limit of $(C_t,\eta_t,p_t,q_t)\in \S_{\Sopt}$ where $C_t$ is smooth. Let $\sigma_t\in \mathbb PH^0(C_t,\eta_t+p_t)$ be a section that vanishes twice at $q_t$ and $\sigma_{X_1},\sigma_{Y_1}$ be aspects of $\sigma_t$ on $X_1$ and $Y_1$ respectively. Then the limit linear series condition tells us that $\ord_{x_1}\sigma_{X_1} + \ord_{y_1}\sigma_{Y_1} \geq g-1$. In the case when $p\in X_1$ and $q\in X_1\cup Y_1$ we find that $q$ is either a double point of $|\eta_{X_1}+p|$ or a double point of $|\eta_{Y_1} + 2y_1|$. In both cases there are finitely many choices for $q$ and the number of possibilities for $q$ is $4g-3$. Recall that the projection $\S_{\Sopt}\to B_1$ is $4g-3$ - degree map, therefore these possibilities exhaust all possible values of $q$ in the case $p\in X_1$. Next, if $p\in Y_1$, $q\in X_1\cup Y_1$ and $k\in \mathbb Z$ is maximal such that $\eta_{Y_1}+p\geq ky_1$ then $q$ is either a double point of $|\eta_{Y_1} + p+y_1|$ or a double point of $|\eta_{X_1} + (k+1)x_1|$. If $k=1$ then there are $4g-1$ possible values for $q$ satisfying this conditions. Excluding $q=p$ and $q=y_1$ we obtain the description of all values of $q$ in the case when $k = 1$. Finally, let $k=2$. Then and $q\notin Y_1$. Then $(X_1\cup_{x_1} E\cup_{y_1} Y_1, \eta^-_{X_1}, \eta_{Y_1}^+, p, q)\in D$ by the definition of $D$. The limit linear series condition applied in this case tells us that either $q\in E$ or $q$ satisfies $\eta_{X_1} + 3x_1\geq 2q$. But this condition trivially holds for each $q\in X_1$. It is easy to show that indeed each $q$ can be realized. We conclude that $\S_{\Sopt}\smm D\to B_1$ is indeed a finite map and $D\subset \S_{\Sopt}$. The fact that $D$ is a divisor in $\S_{\Sopt}$ follows easily.

\end{proof}

Define the divisor $D_{\fib}$ in $\S_{\Sopt}$ by
\begin{equation}
  D_{\fib} := 2D_1 + D_2
  \label{eq:def_of_Dfib}
\end{equation}
where $D_1$ and $D_2$ are divisors from Lemma~\ref{lemma:Dfib_is_fiber}.
\begin{lemma}
  Let $\pi_1: \C_1\to B_1$ be the projection. Then we have
  \begin{align}
    &\S_{\Sopt}\cdot \S_{\dhat} \equiv D_{\W} + D_{\Z}, \label{eq:intersectio_for_DW}\\
    &\S_{\Sopt}\cdot (\S_{\Sopt} + \omega_{\pi_1}) \equiv D_{\Cau} + D_{\W} + D_{\fib}. \label{eq:intersection_for_DCau}
  \end{align}
  \label{lemma:intersection_theory_Gj}
\end{lemma}

\begin{proof}
  The relation~\eqref{eq:intersectio_for_DW} follows from Lemma~\ref{lemma:Zg_cdot}. Let us proof the relation~\eqref{eq:intersection_for_DCau}. Note that the occurrence of the first two summands reflects the adjunction formula applied to $\S_{\Sopt}\subset \C_1$ (recall that $\S_{\Sopt}\to B_1$ has a simple ramification along $D_{\Cau}\cup D_{\W}$ outside the support $D_{\fib}$). More precisely, consider the natural homomorphism $\O_{\C_1}(-\S_{\Sopt})\otimes \omega_{\pi_1}^{\vee}|_{\S_{\Sopt}}\to O_{\S_{\Sopt}}(-D_{\Cau} - D_{\W})$ constructed as follows. Pick a point $[c]\in \S_{\Sopt}$ and let $f\in \O_{\C_1}|_{[c]}$ be a generator of $\O_{\C_1}(-\S_{\Sopt})|_{[c]}$. Let $v\in \omega_{\pi_1}^{\vee}|_{[c]}$ be a generator and represent $v$ as $v_1\otimes v_2$ under the isomorphism $\omega_{\pi_1}\simeq  \pi_1^*K_{B_1}\otimes K_{\C_1}^{\vee}$. Now, send $f\otimes v$ to $(df\wedge v_1)\times v_2\in K_{\C_1}\otimes K_{\C_1}^{\vee}|_{[c]}\simeq \O_{\C_1}|_{[c]}$. It is clear that the image of $f\otimes v$ vanishes along the germ of $D_{\Cau}\cup D_{\W}$ at $[c]$, thus we get a map $\O_{\C_1}(-\S_{\Sopt})\otimes \omega_{\pi_1}^{\vee}|_{\S_{\Sopt}}\to\O_{\C_1}(-D_{\Cau} - D_{\W})$. Since $\S_{\Sopt}\to B_1$ is a ramification cover outside the support of $D_{\fib}$ we conclude that the divisor of the homomorphism is equal to $aD_1 + bD_2$ for some positive integers $a,b$.

  We compute $a$ and $b$ in two steps. First, let us show that $aD_1 + bD_2 = k D_{\fib}$ for some integer $k$. We deduce in from Remark~\ref{rem:stack_to_space_is_ramified}. Namely, let $(X_1\cup_{x_1} E\cup_{y_1} Y_1, \eta^-_{X_1}, \eta_{Y_1}^+, p, q)\in D$ be generic and let $\U\subset \cSo$ be a small neighborhood of the spin curve $(X_1\cup_{x_1} E\cup_{y_1} Y_1, \eta^-_{X_1}, \eta_{Y_1}^+)$ and let $\overline{\C}\to \U$ be the pullback of the universal spin curve from $\cM$ (that is fibers of $\overline{\C}\to \U$ are stable). Denote by $\overline{\C}_u$ the fiber over $u\in \U$ and by $\eta_u$ the corresponding theta characteristics on $\overline{\C}_u$ and let $0\in \U$ be the center of the neighborhood. Note that $(p,y_1)\in T_{\eta_{Y_1}}$. Since $T_{\eta_{Y_1}}$ is smooth and $(p,y_1)$ is a ramification point of the ``second'' projection $T_{\eta_{Y_1}}\to Y_1$ that sends $(y,y')\mapsto y'$ we see that the ``first'' projection $T_{\eta_{Y_1}}\to Y_1$ that sends $(y,y')\mapsto y$ in unbranched at $(p,y_1)$ and thus we can find a continuous family of $p_u\in \supp(\eta_u)$ such that $p_0 = p$. Let
  \begin{equation}
    \overline{\S}_{\Sopt} := \text{closure of } \{(u, q_u)\in \overline{\C}\ \mid\ \overline{\C}_u\text{ is smooth and }\eta_u + p_u\geq 2q_u,\quad q_u\neq p_u\}
    \label{Def_Of_barS}
  \end{equation}
  Then it is clear that the variety $\S_{\Sopt}$ is locally the pullback of $\overline{\S}_{\Sopt}$ and moreover the homomorphism that we constructed above is locally the pullback of the corresponding homomorphism for $\overline{\S}_{\Sopt}$. On the other hand, it follows straighforwardly from the construction of the universal deformation of spin curve that $D_{\fib}$ is locally the pullback of the divisor $\overline D$ of one-dimensional fibers of the projection $\overline{\S}_{\Sopt}\to \U$. It follows that the divisor of the homomorphism constructed above should be a multiple of $D_{\fib}$.

  It remains to show that $k=1$. To do this we will make a local computation. Let $\U_0\subset \U$ be the subvariety of $\U$ parametrizing nodal curves such that for the even component $(Y_u,\eta_{Y_u}^+)$ the nodal point $y_u$ satisfies $\eta_{Y_u}^+ + p_u\geq 2y_u$. Let $V$ be a small polydisc in $\mathbb C^2$ and $\widetilde V$ be the blow up of $V$ at the origin. We construct the map $\widetilde V\times \U_0\to \U$ as follows. Denote by $(X_u\cup_{x_u} E\cup_{y_u} Y_u,\eta_{X_u},\eta_{Y_u})$ the curve corresponding to a point $u\in \U_0$. Let $\sigma_{u1}, \sigma_{u2}\in H^0(X_u, \eta_u + 3x_u)$ be sections such that $\sigma_{u1}$ has a simple zero at $x_u$, $\sigma_{u2}$ does not vanish at $x_u$ and the meromorphic differentials $\sigma^2_{u1}, \sigma^2_{u2}$ corresponding to squares of these sections vary holomorphically with respect to $u$. Moreover, let $\sigma_{u3}\in H^0(Y_u, \eta_{Y_u}+p_u)$ be a section that has a double zero at $y_u$ (such a section is unique up to a constant) chosen such that the corresponding meromorphic differential $\sigma^2_{u3}$ varies holomorphically. Now, consider a point $x:=(t,[\tau_1:\tau_2])\in \widetilde V$. We can find a unique coordinate $z$ at $x_u$ such that $(\tau_1\sigma_{u1}+\tau_2\sigma_{u_2})^2 = z^{-6}(\tau_1z-\tau_2)^2\,dz$ and a unique coordinate $w$ at $y_u$ such that $\sigma^2_{u3} = w^4\,dw$. Moreover, mimic the procedure described in the proof of Lemma~\ref{lemma:tSopt_is_smooth} we can change $(Y_u, y_u)$ locally near $y_u$ to obtain $(Y_{u,v}, y_{u,v})$ and a meromorphc differential $\sigma_{u3,v}^2$ such that $\sigma_{u3,v}^2$ coincides with $\sigma_{u3}^2$ outside a small neighborhood of $y_{u,v}$ and $\sigma_{u_3,v}^2 = -w^2(t\tau_1 - \tau_2w)^2\,dw$ near $y_{u,v}$. Note that if $t = 0$ then $Y_{u,v} = Y_u$. Now consider a two-dimensional family of curves $C_v$ such that $C_v$ is obtained by gluing a neighborhood of $x_u$ and a neighborhood of $y_{u,v}$ via the equation $zw=t$. Note that if $t\neq 0$ then differentials $t^5(\tau_1\sigma_{u1}+\tau_2\sigma_{u_2})^2$ and $\sigma_{u_3,v}^2$ glues to a meromorphic differential $\sigma^2_v$ on $C_v$ with even singularities. Thus we get an odd spin structure $\eta_v$ on $C_v$ and form a map that send $(v,u)\in \widetilde V\times \U_0$ to $(C_v, \eta_v)$. It is clear that the obtained map $\widetilde V\times \U_0\to \U$ is isomorpic to the blow up $\widetilde \U$ of $\U$ along $\U_0$.

  Let $\widetilde \C\to \widetilde \U$ be the pullback of $\overline \C$ and $\varphi: \widetilde \C\to \overline \C$ be the corresponding projection. Set $\widetilde D = \varphi^{-1}(\overline D)$. Then $\varphi^*\overline \S_{\Sopt} = \widetilde \S_{\Sopt} + k\widetilde D$ where $\widetilde \S_{\Sopt}$ is such that the projection $\widetilde \S_{\Sopt}\to \widetilde \U$ is a branched cover. A straighforward analysis of our construction of the blow up made above implies that a pount $q\in X_u$ corresponds to $\widetilde \S_{\Sopt}$ iff it is a ramification point of the one-dimensional linear subsystem of $|\eta_{X_u} + 3x_u|$ spanned by $\tau_1\sigma_{u1}+\tau_2\sigma_{u_2}$ and the section of $\eta_{X_u}$. It follows in particular that the exceptional divisor in $\widetilde \U$ is not an image of the ramification locus of $\widetilde \S_{\Sopt}\to \widetilde \U$ and $k=1$.

  Let $(X_1\cup_{x_1} E\cup_{y_1} Y_1, \eta^-_{X_1}, \eta_{Y_1}^+, p, q)\in D$ be generic. We compute values of $a$ and $b$ constructing a pluming family for this curve.  and let $C_t$ be a plumbing family for the stable curve $X_1\sqcup Y_1/_{x_1\sim y_1}$, i.e. $C_t$ is a smooth one-dimensional family of curves locally given by $\{(z,w,t)\ \mid\ zw = t\}$ such that $C_0$ is isomorphic to $X_1\sqcup Y_1/_{x_1\sim y_1}$ and $C_t$ is smooth for $t\neq 0$. Since the map $\cSo\to \cM$ is unbranched near $(X_1\cup_{x_1} E\cup_{y_1} Y_1, \eta^-_{X_1}, \eta_{Y_1}^+)$ we can find a family of odd theta characteristics $\eta_t$ on $C_t$ such that $\eta_0$ is given by $\eta_{X_1}$ and $\eta_{Y_1}$. Note that points from $\supp(\eta_t)$ converges to $\supp(\eta_{X_1})$ on $X_1$ and to points $y\in Y_1$ such that $(y,y_1)\in T_{\eta_{Y_1}}$. Since $T_{\eta_{Y_1}}$ is smooth and $(p,y_1)$ is a ramification point of the ``second'' projection $T_{\eta_{Y_1}}\to Y_1$ that sends $(y,y')\mapsto y'$ we see that the ``first'' projection $T_{\eta_{Y_1}}\to Y_1$ that sends $(y,y')\mapsto y$ in unbranched at $(p,y_1)$ and thus we can find a continuous family of $p_t\in \supp(\eta_t)$ such that $p_0 = p$. Now let $\sigma_t\in \mathbb PH^0(C_t,\eta_t+p_t)$ be sections with double zeros $q_t$ such that $q_0 = q$ or $q_0 = x_1$ if $q$ belongs to the exceptional component. 
\end{proof}

\begin{proof}[Proof of Lemma~\ref{lemma:Gj_Ddivisor_instead_divisor}]
  As we mentioned before, it is enough to show that $D_{\Cau}$ does not intersect the exceptional component of any fiber of $\pi_1:\C_1\to B_1$ over $Y_i\subset B_1$ for any $i$. Consider first some $i \leq j-1$. Then the claim follows from the analysis of the corresponding limit linear series. For, let $(C, p, q)$, $p,q\in C$, be a quasi-stable marked curve semistably equivalent to a marked spin curve $X\cup_x E\cup_y Y, \eta_X,\O_E(1), \eta_Y, p_{X,i}$ lying in the fiber over $Y_i$. Let $\sigma\in \overline G_g^0(C)$ be an aspect of a limit linear series such that $\sigma_X\in \mathbb PH^0(X, \eta_X\otimes \O_X(p_{X,i} + (g-j)x)), \sigma_X\in \mathbb PH^0(Y, \eta_Y\otimes \O_Y(jy))$ and $\sigma_E\in \mathbb PH^0(E, \O_E(g-1))$. Let moreover $\ord_q \sigma \geq 3$. Notice that the maximal order of $\sigma_X$ at $x$ is $g-j+1$ and the maximal order of $\sigma_Y$ at $y$ is $j-1$. These observations together with the condition that $\sigma$ is a limit linear series implies that $q$ cannot lie on a rational component of $C$. 

  Now assume that $i = j$, so that any fiber of $\pi_1$ over $Y_i$ is isomorphic to a marked spin curve $X\cup_x E\cup_y Y, \eta_X, \O_E(1), \eta_Y, p$ for some $(y,p)\in T_{\eta_Y}$. Let $\iota_{Y, E}: T_{\eta_Y}\times E\hookrightarrow \C_1$ be the underlying embedding of $E$. We will show that $\iota_{Y,E}^*D_{\Cau} = 0$ which implies the statement we are proving. In this purpose we use that
  \begin{equation*}
    \iota_{Y,E}^* D_{\Cau} = \iota_{Y,E}^*(\S_{\Sopt}\cdot (\S_{\Sopt} + \omega_{\pi_1}) - D_{\W} - D_{\fib}).
  \end{equation*}
  Straightforward computations shows that $\iota_{Y,E}^*D_{\fib} = 0$ and $\iota_{Y,E}^*\omega_{\pi_1} = 0$. We also have $\iota_{Y,E}^*D_{\W} = 0$ due to the fact that $D_{\W}$ does not intersect exceptional components. Therefore
  \begin{equation*}
    \iota_{Y,E}^* D_{\Cau} = \iota_{Y,E}^*(\S_{\Sopt}\cdot \S_{\Sopt}).
  \end{equation*}
  We claim that the homology class of $\iota_{Y,E}^*\S_{\Sopt}$ is equal to some constant times the homology class of a vertical fiber of $T_{\eta_Y}\times E$, therefore $\iota_{Y,E}^* D_{\Cau} = \iota_{Y,E}^*(\S_{\Sopt}\cdot \S_{\Sopt}) = 0$. To handle this claim we will show that $\S_{\Sopt}\cap \iota_{Y,E}( \{(y,p)\}\times E) = \varnothing$ for a generic $(y,p)\in T_{\eta_Y}$. We will again use the limit linear series theory. For, let $(C, p, q)$, $p,q\in C$, be a quasi-stable marked curve semistably equivalent to a marked spin curve $X\cup_x E\cup_y Y, \eta_X,\O_E(1), \eta_Y, p_{X,i}$. Let $\sigma\in \overline G_g^0(C)$ be an aspect of a limit linear series such that $\sigma_X\in \mathbb PH^0(X, \eta_X\otimes \O_X(p_{X,i} + (g-j)x)), \sigma_X\in \mathbb PH^0(Y, \eta_Y\otimes \O_Y(jy))$ and $\sigma_E\in \mathbb PH^0(E, \O_E(g-1))$. Let moreover $\ord_q \sigma \geq 2$. We want to show that $q$ does not lies on any rational component of $C$. Notice first that the number (counting multiplicities) of all possible points $q$ is bounded by $\deg(\S_{\Sopt}\to B_1) = 4g-3$. Assume that $q\in X\smallsetminus \{x\}$. Then the limit linear series condition implies that $q$ is a ramification point of the pencil $|\eta_X + 2x|$ not equal to $x$. It follows that there is $4j-1$ such points. Now let $q\in Y\smallsetminus \{y\}$. Limit linear series condition implies that $q$ is a ramification point of $|\eta_Y + y + p|$ different from $p$ and from $y$. Thus Riemann-Hurwitz formula implies that there is $4(g-j)-2$ admissible points $q$. It follows that there is already $4g-3$ possible points $q$ lying away from rational components of $C$ and thus there is no point on rational components.
\end{proof}

\begin{lemma}
  We have $\deg\iota_{X,i,X}^*D_{\Cau} = \deg\iota_{X,i,X}^*D_{\W} = 0$ for each $i = 1,\dots,j$.
  \label{lemma:Gj_XiX}
\end{lemma}

\begin{proof}
  Arguing as in the proof of Lemma~\ref{lemma:Gj_Ddivisor_instead_divisor} we can show that if $(y,q)\in \iota_{X,i,X}^{-1}(D_{\Cau})$ then $\eta_X + p_i \geq 3q$. Since we choose $X$ to be generic we can achieve that there is no such $q\in X$. Similarly, if $(y,q)\in \iota_{X,i,X}^{-1}(D_{\W})$ then $q$ is a base point of $|\eta_X + p_i|$. Again, a generic choice of $X$ ensures that there is no such $q$.
\end{proof}

Now it is left to compute all the intersection numbers in the righ-hand side of~\eqref{eq:DCau_instead_of_Cau_D}~-~\eqref{eq:DZ_instead_of_Z}. Let us start with the following basic lemma. We left the proof of it for the reader.

\begin{lemma}
  Let $(Y,\eta_Y)$ be a generic even spin curve of genus $g-j$ such that the curve $T_{\eta_Y}$ is smooth and connected. Then the projection $T_{\eta_Y}\to Y$ has $4(g-j-1)(g-j)$ ramification points and all of them are simple
  \label{lemma:ramification_of_Scorza}
\end{lemma}

\begin{lemma}
  We have
  \begin{align*}
    \deg\iota_{Y,X}^* D_{\fib} = -4(g-j-1)(g-j)\\
    \deg\iota_{Y,Y}^* D_{\fib} = 4(g-j-1)(g-j)\\
    \deg\iota_{X,i,X}^* D_{\fib} = \deg\iota_{X,i,Y}^* D_{\fib} = 0.
  \end{align*}
  \label{lemma:Gj_Dfib_cap_YX}
\end{lemma}

\begin{proof}
  Follows from Lemma~\ref{lemma:ramification_of_Scorza}.
\end{proof}

\begin{lemma}
  For each $i  = 1,\dots, j-1$ we have
  \begin{align*}
    & \deg\iota_{X,i,Y}^* D_{\W} = 4(g-j-1)(g-j)
  \end{align*}
  \label{lemma:Gj_W_XiY}
\end{lemma}

\begin{proof}
  Consider the correspondence
  \begin{equation}
    \R := \{ (y,q)\in Y\times Y\ \mid\ \eta_Y + 2y \geq 2q \}
    \label{eq:def_of_R}
  \end{equation}
  Studying the corresponding limit linear series we find that $\R = \iota_{X,i, Y}^*\S_{\Sopt}$. Similarly we have $T_{\eta_Y} = \iota_{X,i,Y}^*\S_{\dhat}$. It follows from Lemma~\ref{lemma:intersection_theory_Gj} that
  \begin{equation}
    \iota_{X,i,Y}^* (D_{\W} + D_{\Z}) = \R\cdot T_{\eta_Y}.
    \label{eq:pullback_as_intersection_XiY}
  \end{equation}
  As a straightforward consequence from~\cite[Proposition~5.2]{FARo} we find that
  \begin{equation}
    \deg\iota_{X,i,Y}^*  D_{\Z} = 4(g-j-1)(g-j)
    \label{eq:Gj_DZ_XiY}
  \end{equation}
  Let $E,F\in H_2(Y\times Y)$ denote the classes of a horizontal and a vertical fibers and let $\Delta\in H_2(Y\times Y)$ denote the class of the diagonal. It is well-known that
  \begin{equation}
    T_{\eta_Y} \equiv (g-j-1)(E+F) + \Delta.
    \label{eq:Scorza_in_H2}
  \end{equation}
  Clearly, $\R\cdot \Delta = 0$ and Riemann-Hurwitz formula implies that $\R\cdot E = \R\cdot F = 4(g-j)$. It follows that
  \begin{equation*}
    \deg\R\cdot T_{\eta_Y} = 8(g-j-1)(g-j).
  \end{equation*}
  Substituting this equality and~\eqref{eq:Gj_DZ_XiY} into~\eqref{eq:pullback_as_intersection_XiY} we obtain the desired relation.

\end{proof}

\begin{lemma}
  For each $i  = 1,\dots, j-1$ we have
  \begin{align*}
    & \deg\iota_{X,i,Y}^* D_{\Cau} = 36(g-j-1)(g-j)
  \end{align*}
  \label{lemma:Gj_Cau_XiY}
\end{lemma}

\begin{proof}
  Note that $\iota_{X,i,Y}^*\omega_{\pi_1} \equiv \pr_2^*K_Y + \Delta$ where $\pr_2: Y\times Y\to Y$ is the projection onto the second factor and $\Delta\in H_2(Y\times Y)$ denotes the class of the diagonal. Using this observation and Lemma~\ref{lemma:intersection_theory_Gj} we find that
  \begin{equation}
    \iota_{X,i,Y}^*(D_{\Cau} + D_{\W}) = \R\cdot(\R + \pr_2^*K_Y + \Delta)
    \label{eq:pullback_as_intersection_XiY_Cau+W}
  \end{equation}
  where $\R$ is the correspondence introduced in~\eqref{eq:def_of_R}. Using Seesaw theorem and the fact that $\R$ is symmetrical under the involution on $Y\times Y$ we can write
  \begin{equation}
    \R \equiv d(E + F) + 4\Delta\qquad \text{in } H^2(Y\times Y)
  \end{equation}
  for some $d\in \mathbb Z$ where $E, F, \Delta\in H_2(Y\times Y)$ were denote the classes of a horizontal and a vertical fibers. Riemann-Hurwitz formula implies that $\R\cdot E = 4(g-j)$, therefore
  \begin{equation}
    \R \equiv 4(g-j-1)(E+F) + 4\Delta.
    \label{eq:R_in_H2_XiY}
  \end{equation}
  Using this relation we find that
  \begin{equation*}
    \deg\R\cdot(\R + \pr_2^*K_Y + \Delta) = 40(g-j-1)(g-j).
  \end{equation*}
  Substituting this relation to~\eqref{eq:pullback_as_intersection_XiY_Cau+W} and using that $\deg\iota_{X,i,Y}^*D_{\W} = 4(g-j-1)(g-j)$ due to Lemma~\ref{lemma:Gj_W_XiY} we obtain the desired equality.
\end{proof}

\begin{lemma}
  We have
  \begin{equation*}
    \deg\iota_{Y,X}^*D_{\W} = 4(g-j-1)(j-1)(g-j).
  \end{equation*}
  \label{lemma:Gj_W_YX}
\end{lemma}

\begin{proof}
  Recall that $D_{\W}$ is defined as
  \begin{equation*}
    D_{\W} := \text{closure of } \{(C,\eta, p,q)\in \C_1\ \mid\ C\text{ is smooth and } h^0(C, \eta + p - q) = 2\}.
  \end{equation*}
  The condition $h^0(C,\eta+p-q)=2$ is symmetric in $p,q$, using this observation we will prove that $\deg\iota_{Y,X}^*D_{\W} = \sum_{i = 1}^{j-1}\deg\iota_{X,i,Y}^*D_{\W}$ which leads to the equation in the lemma due to Lemma~\ref{lemma:Gj_W_XiY}.

  Consider $\C_2 := \C\times_B\C$. Then $\C_1$ is a subvariety of $\C_2$ given by
  \begin{equation*}
    \text{closure of } \{(C,\eta,p,q)\in \C_2\ \mid\ C\text{ is smooth and } \eta\geq p\}.
  \end{equation*}
  Note that the variety $D_{\W}$ as a subvarity of $\C_2$ is invariant under the involution $\mathrm{inv}:\C_2\to \C_2$ that interchanges factors. Consider two embeddings
  \begin{align*}
    & \widetilde\iota_{X, Y}: Y\times X\times Y\hookrightarrow \C_2\qquad \widetilde\iota_{X, Y}(y,p,q) = (X\cup_x E\cup_y Y, p, q),\\
    & \widetilde\iota_{Y, X}: Y\times Y\times X\hookrightarrow \C_2\qquad \widetilde\iota_{X,i,X}( (y,p,q) = (X\cup_x E\cup_y Y, p, q).
  \end{align*}
  Note that $\mathrm{inv}\circ \widetilde\iota_{X, Y} = \widetilde\iota_{Y, X}$. On the other hand, we have
  \begin{equation*}
    \sum_{i = 1}^{j-1}\deg\iota_{X,i,Y}^*D_{\W} = \deg\widetilde\iota_{X, Y}^*D_{\W} = \deg\widetilde\iota_{Y, X}^*D_{\W} = \deg\iota_{Y,X}^*D_{\W}
  \end{equation*}
  where the middle equality follows from the fact that $\mathrm{inv}^*D_{\W} = D_{\W}$.
\end{proof}

\begin{lemma}
  We have
  \begin{equation*}
    \deg\iota_{Y,X}^* D_{\Cau} = 4(g-j-1)(g-j)(9j-1).
  \end{equation*}
  \label{lemma:Gj_Cau_YX}
\end{lemma}

\begin{proof}
  Introduce two subvarieties $\Efib, \Ffib\subset T_{\eta_Y}\times X$ given by
  \begin{align*}
    \Efib := \{  ( (y,p), q)\in T_{\eta_Y}\times X\ \mid\ \eta_X + 2x \geq 2q,\ q\neq x \},\\
    \Ffib := \{ ( (y,p), q)\in T_{\eta_Y}\times X\ \mid\ \eta_Y + p\geq 2y \}
  \end{align*}
  (recall that $x\in X$ is the nodal point). It follows from Lemma~\ref{lemma:Dfib_is_fiber} that
  \begin{equation*}
    \iota_{Y,X}^*\S_{\Sopt} = \Efib + \Ffib
  \end{equation*}
  Notice that $\Efib\equiv (4j-1)E$, $\Ffib \equiv 4(g-j-1)(g-j)F$ where $E,F\in H_2(T_{\eta_Y}\times X)$ are classes of a horizontal and a vertical fibers. We also have $\iota_{Y,X}^*\omega_{\pi_1} \equiv \pr_2^*K_X + E$ where $\pr_1: T_{\eta_Y} \times X\to T_{\eta_Y} $, $\pr_2: T_{\eta_Y} \times X\to X$ are projections onto the first and the second factors. It follows that
  \begin{equation*}
    \begin{split}
      &\deg(\iota_{Y,X}^*(D_{\Cau} + D_{\W})) -4(g-j-1)(g-j) = \deg\iota_{Y,X}^* (\S_{\Sopt}\cdot(\S_{\Sopt} + \omega_{\pi_1})) =\\
      &= ((4j-1)E + 4(g-j-1)(g-j)F)\cdot ((6j-2)E + 4(g-j-1)(g-j)F)=\\
      &= 4(g-j-1)(g-j)(10j-3).
    \end{split}
  \end{equation*}
  thus
  \begin{equation*}
    \deg(\iota_{Y,X}^*D_{\Cau} = 4(g-j-1)(g-j)(10j-3) - \deg(\iota_{Y,X}^*D_{\W}.
  \end{equation*}
  Using Lemma~\ref{lemma:Gj_W_YX} we conclude that $\deg\iota_{Y,X}^* D_{\Cau} = 4(g-j-1)(g-j)(9j-1)$.
\end{proof}

\begin{lemma}
  We have 
  \begin{equation*}
    \deg\iota_{Y,Y}^*D_{\W} = 2(g-j-1)(g-j)(g-j+1).
  \end{equation*}
  \label{lemma:Gj_DW_YY}
\end{lemma}

\begin{proof}
  Introduce the correspondence $\R\subset T_{\eta_Y}\times Y$ given by
  \begin{equation*}
    \begin{split}
      \R := & \text{ closure of }\{  ((y, p), q)\in T_{\eta_Y}\times Y\ \mid\ \eta_Y + p + y\geq 2q,\ p\neq q,\ y\neq q \}\\
      = & \{  ((y,p), q)\in T_{\eta_Y}\times Y\ \mid\ \eta_Y + p + y\geq 2q,\ p\neq q,\ y\neq q \}\cup \\
      & \cup \{  ((p, q), q)\in T_{\eta_Y}\times Y\ \mid\ \eta_Y + p \geq 2q,\} \cup \{  ((y, q), q)\in T_{\eta_Y}\times Y\ \mid\ \eta_Y + y \geq 2q,\}
    \end{split}
  \end{equation*}
  It is straightforward that $\iota^*_{Y,Y}\S_{\Sopt} = \R$. Let $E,F,\Delta_y,\Delta_p\in H_2(T_{\eta_Y}\times Y, \mathbb Z)$ be the classes of horizontal and vertical fibers and the diagonals, i.e.
  \begin{equation*}
    \begin{split}
      & E = \{ ((y,p), q_0)\in T_{\eta_Y}\times Y\ \mid\ q_0\text{ is fixed} \},\\
      & F = \{ ( (y_0,p_0), q)\in T_{\eta_Y}\times Y\ \mid\ (y_0,p_0)\text{ is fixed} \},\\
      & \Delta_y = \{ ((q,y), y)\ \mid\ (q, y)\in T_{\eta_Y} \},\\
      & \Delta_p = \{ ((y,p), p)\ \mid\ (y, p)\in T_{\eta_Y} \}.
    \end{split}
  \end{equation*}
  Using Riemann-Hurwitz formula we find
  \begin{equation}
    \R\cdot F = 4(g-j) - 2.
    \label{eq:Gj_R_cdot_F_YY}
  \end{equation}
  To compute the intersection number $\R\cdot E$ we embedd $\R$ into $Y^3$ under $T_{\eta_Y}\times Y\hookrightarrow Y^3$. Consider the following varieties
  \begin{equation*}
    \begin{split}
      & Q_{T_{\eta_Y}} = \{(y, p, q)\in Y^3\ \mid\ (y, p) \in T_{\eta_Y} \},\\
      & Q_2 = \{(y, p, q)\in Y^3\ \mid\ T_{\eta_Y} + y + p \geq 2q \},\\
      & V_0 = \{(y, p, q_0)\in Y^3\ \mid\ q_0\text{ is fixed} \}.
    \end{split}
  \end{equation*}
  It is straightforward that $(\R+\Delta_y + \Delta_p)\cdot E= Q_{T_{\eta_Y}}\cdot Q_2 \cdot V_0$. Now let $Q\subset Y^2$ be defined
  \begin{equation*}
    Q = \{(y, p)\in Y^2\ \mid\ \eta_Y + y + p\geq 2q_0\}.
  \end{equation*}
  We have
  \begin{equation*}
    \deg(Q_{T_{\eta_Y}}\cdot Q_2 \cdot V_0) = \deg(Q\cdot T_{\eta_Y}).
  \end{equation*}
  Using that the relation~\eqref{eq:Scorza_in_H2} we find that
  \begin{equation*}
    Q\cdot T_{\eta_C} = 2(g-j)(g-j+1)
  \end{equation*} 
  and conclude that
  \begin{equation}
    \R\cdot E = 2(g-j)^2.
    \label{eq:Gj_R_cdot_E_YY}
  \end{equation}
  A direct application of Seesaw Theorem gives that 
  \begin{equation}
    \R \equiv a E + b F + \Delta_y + \Delta_p
    \label{eq:Gj_Seesaw_for_R_YY}
  \end{equation}
  for some integers $a,b$. Combaining~\eqref{eq:Gj_R_cdot_F_YY},~\eqref{eq:Gj_R_cdot_E_YY} and~\eqref{eq:Gj_Seesaw_for_R_YY} we conclude that
  \begin{equation}
    \R \equiv 4(g-j-1)\, E + 2(g-j)(g-j-1)\, F + \Delta_y + \Delta_p.
    \label{eq:R_in_H2_YY}
  \end{equation}
  
  Let $\pr_y: T_{\eta_Y}\times Y \to Y\times Y$ be the projection defined by $\pr_y( (y,p), q)  = (y,q)$. We have $\iota^*_{Y,Y}\S_{\dhat} = \pr_y^*T_{\eta_Y} - \Delta_p$. It follows from Lemma~\ref{lemma:intersection_theory_Gj} that
  \begin{equation*}
    \deg\iota_{Y,Y}^* D_{\W} = (\pr_y^*T_{\eta_Y} - \Delta_p)\cdot \R - \iota_{Y,Y}^*D_{\Z}.
  \end{equation*}
  Using~\eqref{eq:R_in_H2_YY} we find that $(\pr_y^*T_{\eta_Y} - \Delta_p)\cdot \R = 2(g-j-1)(g-j)(g-j+1) + 4(g-j-1)^2(g-j)$. On the other hand, we have $\iota_{Y,Y}^*D_{\Z} = 4(g-j-1)^2(g-j)$. It follows that $\deg\iota_{Y,Y}^* D_{\W} = 2(g-j-1)(g-j)(g-j+1)$ as required.
\end{proof}

\begin{lemma}
  We have
  \begin{equation*}
    \deg\iota_{Y,Y}^*D_{\Cau} = 2(g-j-1)(g-j)(9(g-j)+1).
  \end{equation*}
  \label{lemma:Gj_DCau_YY}
\end{lemma}

\begin{proof}
  Let us follow the notation of Lemma~\ref{lemma:Gj_DW_YY}. Using Lemma~\ref{lemma:intersection_theory_Gj} and the relation $\iota^*_{Y,Y}\S_{\Sopt} = \R$ we obtain
  \begin{equation*}
    \deg\iota_{Y,Y}^*D_{\Cau} = \R\cdot (\R + \iota_{Y,Y}^*\omega_{\pi_1}) - \deg\iota_{Y,Y}^*(D_{\W} + D_{\fib}).
  \end{equation*}
  Recall that $E\in H_2(T_{\eta_Y}\times Y)$ denotes the class of a horizontal fiber. We have $\iota_{Y,Y}*\omega_{\pi_1}\equiv (2(g-j)-2)E + \Delta_y$ in $H^2(T_{\eta_Y}\times Y)$. Using this observation and~\eqref{eq:R_in_H2_YY} we find that
  \begin{equation*}
    \R + \iota_{Y,Y}^*\omega_{\pi_1} \equiv 6(g-j-1)\, E + 2(g-j)(g-j-1)\, F + 2\Delta_y + \Delta_p.
  \end{equation*}
  Using this relation and~\eqref{eq:R_in_H2_YY} again we find that
  \begin{equation*}
    \R\cdot (\R + \iota_{Y,Y}^*\omega_{\pi_1}) = 2(g-j-1)(g-j)(10(g-j) + 4).
  \end{equation*}
  Recall that by Lemma~\ref{lemma:Gj_Dfib_cap_YX}
  \begin{equation*}
    \deg\iota_{Y,Y}^*D_{\fib} = 4(g-j-1)(g-j).
  \end{equation*}
  Summing up computations above and using the result of Lemma~\ref{lemma:Gj_DW_YY} we find that
  \begin{equation*}
    \deg\iota_{Y,Y}^*D_{\Cau} = 2(g-j-1)(g-j)(9(g-j)+1).
  \end{equation*}
\end{proof}

\begin{rem}
  Let us comment a little bit on results of Lemma~\ref{lemma:Gj_W_XiY} - Lemma~\ref{lemma:Gj_DCau_YY}. Denote the curve $(X\cup_x E\cup_y Y, \eta_X, \eta_Y)$ by $(C(y),\eta)$ for simplicity. A study of the corresponding limit linear series can be used to show that
  \begin{equation*}
    G_j\cap \cCau = \CMcal{A}_{XY}\cup \CMcal{A}_{YX}\cup \CMcal{A}_{YY},
  \end{equation*}
  where
  \begin{align}
    \CMcal{A}_{XY} = & \bigcup_{p\in \supp(\eta_X)} \{(C(y), \eta)\ \mid\ \exists q\in Y:\ \eta_Y + 2y \geq 3q \}, \label{eq:XY}\\
    \CMcal{A}_{YX} = & \bigcup_{\substack{ q\in X:\ \eta_X + 3x \geq 3q\\ q\neq x }}\{ (C(y), \eta)\ \mid\ \exists p\in \supp(|\eta_Y + y|):\ \eta_Y + p\geq 2y \}, \label{eq:YX}\\
    \CMcal{A}_{YY} = & \{(C(y), \eta)\ \mid\ \exists p\in \supp(|\eta_Y + y|), q\in Y:\ \eta_Y + p + y\geq 3q,\ q\notin\{y,p\}\} \label{eq:YY}.
  \end{align}
  Unions in the right-hand side of~\eqref{eq:XY} and~\eqref{eq:YX} are superfluous but we write them to emphasize that the multiplicity of the intersection is at least the size of the index set of the union. Then it is possible to show that the number from Lemma~\ref{lemma:Gj_Cau_XiY} is equal to the size of the set~\eqref{eq:XY}, the number from Lemma~\ref{lemma:Gj_Cau_YX} is equal to the size of the set~\eqref{eq:YX} and the number from Lemma~\ref{lemma:Gj_DCau_YY} is equal to the size of the set~\eqref{eq:YY}.

  The same remark can be made in the case of $\cW$.
  \label{rem:interpretaion_via_limit_series}
\end{rem}

\subsection{Intersection of $H_2,\dots, H_{g-2}$ with caustic and fixed point divisor.}\label{subsec:H_cap_Cau_and_W}

Fix some $j,\ 2\leq j\leq g-2$, till the end of this subsection. We compute the intersection number $H_j\cdot \cCau$ and $H_j\cdot \cW$ in a similar way we passed through in Section~\ref{subsec:G_cap_Cau_and_W}. The family $H_j$ is the image of the map of the elliptic curve $T\to \cSo$ that sends $t_2\in T$ to the moduli of the curve $(X\cup_x E_X\cup_{t_1} T\cup_{t_2} E_Y\cup_y Y, \eta_X, \eta_T, \eta_Y)$. Let $\mathbb D\subset \mathbb C^{3g-4}$ be a polydisc and set $B = T\times \mathbb D$. Then we can extend the map $T\to \cSo$ to a map $\CMcal{M}: B\to \cSo $ such that for each $t_2\in T$ there exists a neighborhood $U\subset T$ of $t_2$ such that the pullback of the universal spin curve under the induced map $U\times \mathbb D\to \cSo$ is a universal deformation of $(X\cup_x E\cup_y Y, \eta_X, \eta_Y)$. Denote by $l : \C\to B$ the pullback of the universal spin curve under $\CMcal M$.

Denote by $B_{\cSop}\subset \C$ the closure of the preimage of $\tSop$ under the map $\C\to \cCo$. Shrinking $\mathbb D$ if necessary we can achieve that the variety $B_{\cSop}$ is smooth and has $g-1$ connected components. To describe these components assume that $\eta_X\otimes\O_X(x) \simeq \O_X(p_1 + \dots + p_{j-1})$ and $\eta_Y\otimes\O_Y(y) \simeq \O_Y(p_j + \dots + p_{g-1})$. Then $l^{-1}(Y\times \{0\})\cap B_{\cSop} = T_1\sqcup \dots\sqcup T_{g-1}$ where
\begin{equation*}
  T_i = \{ (X\cup_x E_X\cup_{t_1} T \cup_{t_2} E_Y \cup_y Y, \eta_X, \eta_t, \eta_Y, p_i)\in \C \}_{y\in Y}, \qquad i = 1,\dots, g-1.
\end{equation*}
Denote by $B_i$ the connected component of $B_{\cSop}$ that contains $T_i$ for each $i = 1,\dots, g-1$. If $\mathbb D$ is small enough then $B_{\cSop} = B_1\sqcup\dots\sqcup B_{g-1}$ and the projection $B_i\to B$ is an isomorphism for each $i$.

Let $\C_{\cSop} = \C\times_B B_{\cSop}$ and consider the following codimension $2$ locuses of $\C_{\cSop}$:
\begin{align}
  & D_{\Cau} := \text{closure of } \{(C,\eta, p, q)\in \C_{\cSop}\ \mid\ C\text{ is smooth and } \eta + p \geq 3q,\ q\neq p\},\\
  & D_{\W} := \text{closure of } \{(C,\eta, p,q)\in \C_{\cSop}\ \mid\ C\text{ is smooth and } h^0(C, \eta + p - q) = 2\},\\
  & D_{\Z} := \text{closure of } \{(C,\eta, p,q)\in \C_{\cSop}\ \mid\ C\text{ is smooth and } \eta \geq p+2q\}.
  \label{eq:Hj_def_of_DSubscript}
\end{align}
Let $\C_i\to B_i$ be the restriction of $\C_{\cSop}\to B_{\cSop}$ to $B_i$. Introduce the following embeddings:
\begin{align*}
  & \iota_{i, X}: T\times X\hookrightarrow \C_i\qquad \iota_{i,X}(t,q) = (X\cup_x E_X\cup_{t_1} T\cup_{t_2} E_Y\cup_y Y, p_i, q),\\
  & \iota_{i, Y}: T\times Y\hookrightarrow \C_i\qquad \iota_{X,i,X}(t,q) = (X\cup_x E_X\cup_{t_1} T\cup_{t_2} E_Y\cup_y Y, p_i, q),\\
\end{align*}

\begin{lemma}
  We have
  \begin{align}
    & \cCau \cdot H_j = \sum_{i = 1}^{j-1} \deg\iota_{i,Y}^*D_{\Cau} + \sum_{i = j}^{g-1}\deg\iota_{i,X}^*D_{\Cau},  \label{eq:Hj_DCau_instead_of_Cau_D}\\
    & \cW \cdot G_j = \frac{1}{2}\left( \sum_{i = 1}^{j-1} \deg\iota_{i,Y}^*D_{\W} + \sum_{i = j}^{g-1}\deg\iota_{i,X}^*D_{\W} \right),  \label{eq:Hj_DW_instead_of_W}\\
    & (g-3)\cZ \cdot G_j = \sum_{i = 1}^{j-1} \deg\iota_{i,Y}^*D_{\Z} + \sum_{i = j}^{g-1}\deg\iota_{i,X}^*D_{\Z}.  \label{eq:Hj_DZ_instead_of_Z}
  \end{align}

  \label{lemma:Hj_Ddivisor_instead_divisor}
\end{lemma}
The proof of Lemma~\ref{lemma:Hj_Ddivisor_instead_divisor} is similar to the proof of Lemma~\ref{lemma:Gj_Ddivisor_instead_divisor}, so we omit it. Set
\begin{align}
  & \S_{\Sopt} := \text{closure of } \{(C,\eta, p, q)\in \C_{\cSop}\ \mid\ C\text{ is smooth and } \eta + p \geq 2q,\ q\neq p\},\\
  & \S_{\dhat} := \text{closure of } \{(C,\eta, p, q)\in \C_{\cSop}\ \mid\ C\text{ is smooth and } \eta \geq p+q\}.
  \label{eq:Hj_def_of_SSubscript}
\end{align}
Notice that the variety $\S_{\Sopt}$ corresponds to the variety $\pi_2(\tSopt)$ (see Section~\ref{subsec:relation_between_Cau_and_lambda} and in particular Corollary~\ref{cor:cor_tSopt_is_smooth}). The forgetful map $\S_{\Sopt}\to B_{\cSop}$ is a branched cover outside a locus in the boundary of $B_{\cSop}$ that is described below:

\begin{lemma}
  For each $i = 1,\dots,j-1$ introduce the locus $D_{X,i}\subset \C_i$ given by
  \begin{align*}
    D_X := \{ (X_0\cup_{x_0} E_{X_0}\cup_{t_1} T_0\cup_{t_2} E_{Y_0}\cup_{y_0} Y_0, \eta^+_{X_0}, \eta^-_{T_0}, \eta_{Y_0}^+, p_i, q)\in \C_i\ \mid\\
    \eta_{T_0} + 2t_1 \geq 2t_2,\ t_1\neq t_2,\ q\in Y_0\cup E_{Y_0} \}.
  \end{align*}
  For each $i = j,\dots,g-1$ consider the locus $D_{Y,i}\subset \C_i$ given by
  \begin{align*}
    D_Y := \{ (X_0\cup_{x_0} E_{X_0}\cup_{t_1} T_0\cup_{t_2} E_{Y_0}\cup_{y_0} Y_0, \eta^+_{X_0}, \eta^-_{T_0}, \eta_{Y_0}^+, p, q)\ \mid\\
    \eta_{T_0} + 2t_1 \geq 2t_2,\ t_1\neq t_2,\ q\in X_0\cup E_{X_0} \}.
  \end{align*}
  Set
  \begin{equation*}
    D := \left(\bigcup_{i = 1}^{j-1}D_{X,i} \right) \cup \left( \bigcup_{i = j}^{g-1} D_{Y,i} \right).
  \end{equation*}
  Then $D$ is a divisor in $\S_{\Sopt}$ and the forgetful map $\S_{\Sopt}\to B_1$ is a branched cover outside $D$. As a divisor $D$ splits into a sum of two components $D_1+D_2$ where $D_1$ corresponds to the case when $q\in X_0$ or $q\in Y_0$ while $D_2$ corresponds to the case when $q\in E_F$ or $q\in E_R$.
  \label{lemma:Hj_Dfib_is_fiber}
\end{lemma}

The proof of Lemma~\ref{lemma:Hj_Dfib_is_fiber} is similar to the proof of Lemma~\ref{lemma:Dfib_is_fiber} and we omit it.

Using the notation of Lemma~\ref{lemma:Hj_Dfib_is_fiber} define the divisor $D_{\fib}\subset \S_{\Sopt}$ by
\begin{equation}
  D_{\fib} := 2D_1 + D_2
  \label{eq:Hj_def_of_Dfib}
\end{equation}
\begin{lemma}
  Let $\pi_1: \C_{\cSop}\to B_{\cSop}$ be the projection. Then we have
  \begin{align}
    &\S_{\Sopt}\cdot \S_{\dhat} = D_{\W} + D_{\Z}, \label{eq:Hj_intersectio_for_DW}\\
    &\S_{\Sopt}\cdot (\S_{\Sopt} + \omega_{\pi_1}) = D_{\Cau} + D_{\W} + D_{\fib}. \label{eq:Hj_intersection_for_DCau}
  \end{align}
  \label{lemma:intersection_theory_Hj}
\end{lemma}

The proof of Lemma~\ref{lemma:intersection_theory_Hj} is similar to the proof of Lemma~\ref{lemma:intersection_theory_Gj}, so we omit it.

\begin{lemma}
  For any $i = 1,\dots, j-1$ we have
  \begin{align*}
    &\deg\iota_{i, Y}^*D_{\W} = 3(g-j),\\
    &\deg\iota_{i,X}^*D_{\W} = 0.
  \end{align*}
  For any $i = j,\dots, g-1$ we have
  \begin{align*}
    &\deg\iota_{i, X}^*D_{\W} = 3(j-1),\\
    &\deg\iota_{i,Y}^*D_{\W} = 0.
  \end{align*}

  \label{lemma:Hj_cdot_DW}
\end{lemma}

\begin{proof}
  It is enough to prove the lemma for $i \leq j-1$ and all possible values of $j$: then we can handle the case $i \geq j$ by symmetry: one needs to replace $H_j$ with $H_{g-j+1}$ that is homotopic to $H_j$. Thus let us assume that $i \leq j-1$. The fact that $\deg\iota_{i,X}^*D_{\W} = 0$ follows from $\iota_{i,X}(T\times X)\cap D_{\W}=\varnothing$. Indeed, analyzing the limit linear series condition as in the proof of Lemma~\ref{lemma:Gj_Ddivisor_instead_divisor} we find that if $\iota_{i,X}(t,q)\in D_{\W}$ then $q$ is a base point of the linear system $|\eta_X + x + p_i|$. But as $x$ was chosen generic we may assume $|\eta_X + x + p_i|$ does not have base points.
  
  By the same reasons we get that $\deg\iota_{i,Y}^*D_{\Z} = 0$, thus 
  \begin{equation*}
    \deg\iota_{i, Y}^*D_{\W} = \deg\iota_{i, Y}^*(\S_{\Sopt}\cdot \S_{\dhat})
  \end{equation*}
  due to~\eqref{eq:Hj_intersectio_for_DW}. Let $E,F\in H_2(T\times Y)$ denote the classes of a horizontal and a vertical fibers respectively. It is straightforward that
  \begin{equation*}
    \iota_{i,Y}^* \S_{\dhat}= \{ (t,q)\in T\times Y\ \mid\ \eta_Y + y \geq q \}
  \end{equation*}
  thus
  \begin{equation*}
    \iota_{i,Y}^* \S_{\dhat} \equiv (g-j)E
  \end{equation*}
  in $H_2(T\times Y)$. On the other hand we have
  \begin{equation*}
    \iota_{i,Y}^*\S_{\Sopt} \equiv \{(t,q)\in T\times Y\ \mid\ \eta_Y + 2y\geq 2q\} + \{(t,q)\in T\times Y\ \mid\ \eta_T + 2t_1\geq 2t_2,\ t_1\neq t_2\}
  \end{equation*}
  which implies that
  \begin{equation}
    \iota_{i,Y}^*\S_{\Sopt} \equiv 4(g-j)E + 3F.
    \label{eq:Hj_SSopt_in_H2}
  \end{equation}
  \begin{equation*}
    \deg\iota_{i, Y}^*D_{\W} = (g-j)E\cdot (4(g-j) E + 3F) = 3(g-j).
  \end{equation*}
\end{proof}

\begin{lemma}
  For any $i = 1,\dots, j-1$ we have
  \begin{align*}
    &\deg\iota_{i, Y}^*D_{\Cau} = 27(g-j),\\
    &\deg\iota_{i,X}^*D_{\Cau} = 0.
  \end{align*}
  For any $i = j,\dots, g-1$ we have
  \begin{align*}
    &\deg\iota_{i, X}^*D_{\Cau} = 27(j-1),\\
    &\deg\iota_{i,Y}^*D_{\Cau} = 0.
  \end{align*}

  \label{lemma:Hj_cdot_DCau}
\end{lemma}

\begin{proof}
  Again, it is enough to prove the lemma for $i \leq j-1$ and all possible values of $j$. The fact that $\deg\iota_{i,X}^*D_{\Cau} = 0$ follows from $\iota_{i,X}(T\times X)\cap D_{\Cau} = \varnothing$. Indeed, as in the proof of Lemma~\ref{lemma:Hj_cdot_DW} we can deduce that if $\iota_{i,X}(t,q)\in D_{\Cau}$ then $\eta_X + x + p_i\geq 3q$. But as $x$ was chosen generic we may assume there is no such $q\in X$.
  
  Due to~\eqref{eq:Hj_intersection_for_DCau} we have
  \begin{equation*}
    \deg\iota_{i, Y}^*D_{\Cau} = \deg\iota_{i,Y}*(\S_{\Sopt}\cdot(\S_{\Sopt} + \omega_{\pi_1})) - \deg\iota_{i,Y}^*D_{\W} - \deg\iota_{i,Y}^*D_{\fib}.
  \end{equation*}
  It follows from the construction of $D_{\fib}$ (cf.~Lemma~\ref{lemma:Hj_Dfib_is_fiber}) that $\deg\iota_{i,Y}^*D_{\fib} = -3$. Let $E,F\in H_2(T\times Y)$ denote the classes of a horizontal and a vertical fibers respectively. We have $\iota_{i,Y}*\omega_{\pi_1} \equiv (2(g-j)-1)E$. Thus, using~\eqref{eq:Hj_SSopt_in_H2} and the result of Lemma~\ref{lemma:Hj_cdot_DW} we obtain
  \begin{equation*}
    \deg\iota_{i, Y}^*D_{\Cau} = (4(g-j) E + 3F)\cdot( (6(g-j)-1) E + 3F) - 3(g-j) + 3 = 27(g-j).
  \end{equation*}
\end{proof}

\subsection{Values of the coefficients at $\alpha_j$ for non-zero $j$.} 
Recall the notation for coefficients~\eqref{eq:notation_for_coefficietns}:
\begin{equation*}
  \begin{split}
    & \cCau \equiv l^c \cdot \lambda - a_0^c \cdot \alpha_0 - b_0^c\cdot \beta_0 - \sum_{j = 1}^{g-1}a_j^c \cdot \alpha_j,\\
    & \cW \equiv l^w \cdot \lambda - a_0^w \cdot \alpha_0 - b_0^w\cdot \beta_0 - \sum_{j = 1}^{g-1}a_j^w \cdot \alpha_j\\
  \end{split}
\end{equation*}

\begin{prop}
  For any $j = 1,\dots, g-1$ the following formulas holds:
  \begin{align*}
    &a_j^c = (g-j)(9g + 27j - 19)\\
    &a_j^w = \frac{1}{2}\cdot(g-j)(g + 3j - 3).
  \end{align*}
  \label{prop:alphaj_for_j>0}
\end{prop}

\begin{proof}
  First, suppose that $j\leq g-2$. It follows from~\eqref{eq:Gj_cdot_basis} that
  \begin{equation*}
    2(g-j-1)a_j^c = G_j\cdot \cCau.
  \end{equation*}
  Combining this with~\eqref{eq:DCau_instead_of_Cau_D} we find that
  \begin{align*}
    & 2(g-j-1)a_j^c =\\
    &= \sum_{i = 1}^j (\deg\iota_{X,i,X}^*D_{\Cau} + \deg\iota_{X,i,Y}^*D_{\Cau}) + \deg\iota_{Y,X}^*D_{\Cau} + \deg\iota_{Y,Y}^*D_{\Cau}.
  \end{align*}
  Applying Lemma~\ref{lemma:Gj_XiX}, Lemma~\ref{lemma:Gj_Cau_XiY}, Lemma~\ref{lemma:Gj_Cau_YX} and Lemma~\ref{lemma:Gj_DCau_YY} we compute the right-hand side of this equation. Similarly, due to~\eqref{eq:DW_instead_of_W} we have
  \begin{align*}
    & 2(g-j-1)a_j^w =\\
    &= \frac{1}{2}\left( \sum_{i = 1}^j (\deg\iota_{X,i,X}^*D_{\W} + \deg\iota_{X,i,Y}^*D_{\W}) + \deg\iota_{Y,X}^*D_{\W} + \deg\iota_{Y,Y}^*D_{\W} \right).
  \end{align*}
  Applying Lemma~\ref{lemma:Gj_XiX}, Lemma~\ref{lemma:Gj_W_XiY}, Lemma~\ref{lemma:Gj_W_YX} and Lemma~\ref{lemma:Gj_DW_YY} we obtain the desired formula for $a_j^w$.

  Now, let $j = g-1$. Due to~\eqref{eq:intersection_with_H} we have
  \begin{equation*}
    \begin{split}
      & H_{g-1}\cdot \lambda = H_{g-1} \cdot \beta_0 = 0,\qquad H_{g-1} \cdot\alpha_i = 0,\quad i \neq 1,2,g-1\\
      & H_{g-1}\cdot \alpha_{g-1} = H_{g-1}\cdot \alpha_2 = -1,\\
      & H_{g-1}\cdot \alpha_1 = 1.
    \end{split}
  \end{equation*}
  Using Lemma~\ref{lemma:Hj_Ddivisor_instead_divisor} and Lemma~\ref{lemma:Hj_cdot_DW}, Lemma~\ref{lemma:Hj_cdot_DCau} we find that
  \begin{align*}
    & a_{g-1}^w + a_2^w - a_1^w = 6(g-2),\\
    & a_{g-1}^c + a_2^c - a_1^c = 54(g-2).
  \end{align*}
  Substituting values of $a_1^w, a_2^w, a_1^c$ and $a_2^c$ that we computed above we obtain relations defining $a_{g-1}^w$ and $a_{g-1}^c$.
\end{proof}

\begin{rem}
  It follows from~\eqref{eq:intersection_with_H} that 
  \begin{align*}
    & a_j^w + a_{g-j+1}^w - a_1^w = \cW\cdot H_j,\\
    & a_j^c + a_{g-j+1}^c - a_1^c = \cCau\cdot H_j.
  \end{align*}
  for any $j=2,\dots, g-2$. Using Lemma~\ref{lemma:Hj_Ddivisor_instead_divisor}, Lemma~\ref{lemma:Hj_cdot_DW} and Lemma~\ref{lemma:Hj_cdot_DCau} we conclude that:
  \begin{align*}
    & a_j^w + a_{g-j+1}^w - a_1^w = 6(j-1)(g-j),\\
    & a_j^c + a_{g-j+1}^c - a_1^c = 54(j-1)(g-j).
  \end{align*}
  This is consistent with Proposition~\ref{prop:alphaj_for_j>0}, so it agrees with the formulas for the Caustic and Base Point divisors that we obtained.
  \label{rem:test_via_intersection}
\end{rem}

\subsection{Coefficients at $\alpha_0$ and $\beta_0$.}\label{subsec:alpha0_beta0}

At this point we have computed all coefficients in~\eqref{eq:notation_for_coefficietns} except those at $\alpha_0$ and $\beta_0$. We will compute these two coefficients via intersection numbers of $\cW$ and $\cCau$ with two elliptic pencils $F_0$ and $G_0$ exactly as is was done in~\cite[Proposition~5.3]{FARo} in the case of $\Z$. Let us briefly recall the construction of $F_0$ and $G_0$. Consider the pencil of plane cubic curves $\{ Y_{\lambda} = f^{-1}(\lambda) \}_{\lambda\in \mathbb P^1}$ where $f: \Bl_9(\mathbb P^2)\to \mathbb P^1$. Let $y: \mathbb P^1\to \Bl_9(\mathbb P^2)$ be a section given by one of base points of the pencil of cubics. Let $(X,x)\in \CMcal{M}_{g-1,1}$ be a generic pointed curve and let $\eta_X^+$ and $\eta_X^-$ be even and odd theta-characteristics on it. Define
\begin{equation}
  F_0 := \{X\cup_x E\cup_{y(\lambda)} Y, \eta^+_X, \O_E(1), \eta_{f^{-1}(\lambda)} = \O_{f^{-1}(\lambda)}\} \subset \cSo.
  \label{eq:def_of_F0}
\end{equation}
We have $F_0\cdot \lambda = 1, F_0\cdot \alpha_0 = 12$ and $F_0\cdot \alpha_1 = -1$. Moreover, $F_0\cdot \beta_0 = F_0\cdot \alpha_j = 0$ for $j\geq 2$. Then, define $G_0$ as
\begin{equation}
  G_0 := \{X\cup_x E\cup_{y(\lambda)} Y, \eta^-_X, \O_E(1), \eta_{f^{-1}(\lambda)}\text{ is even} \} \subset \cSo.
  \label{eq:def_of_G0}
\end{equation}
We have $G_0\cdot \lambda = 3$ and $G_0\cdot \alpha_{g-1} = -3$, $G_0\cdot \alpha_0 = G_0\cdot \beta_0 = 12$. Moreover, we have $G_0\cdot \alpha_j = 0$ for each $j = 1,\dots, g-2$.

\begin{lemma}
  \begin{align*}
    & \cW\cdot F_0 = \cW\cdot G_0 = 0,\\
    & \cCau\cdot F_0 = \cCau \cdot G_0 = 0.
  \end{align*}
  \label{lemma:Cau_and_W_cdot_F0_and_G0}
\end{lemma}

\begin{proof}
  The proof easily follows from the analysis of the corresponding limit linear series. Let us prove that $\cW\cdot G_0 = 0$ as an example. For, assume that $(X\cup_x E\cup_{y(\lambda)} Y, \eta^-_X, \O_E(1), \eta_{f^{-1}(\lambda)})\in G_0\cap \cW$. Let $(C, p, q)$, $p,q\in C$, be a quasi-stable marked curve semistably equivalent to marked spin curve $X\cup_x E\cup_{y(\lambda)} Y$. Let $\sigma\in \overline G_g^0(C)$ be an aspect of a limit linear series such that $\sigma_X\in \mathbb PH^0(X, \eta_X\otimes \O_X(x))$ and $(p,q)$ be limits of points $p',q'$ on smooth odd spin curve $(C', \eta')$ such that $h^0(C', \eta' + p' - q') = 2$. Note that the condition on $p'$ and $q'$ is symmetrical. At least one of the points $p,q$ lies on $X$, so we can assume that $q\in X$ without loss of generality. Then the condition on $\sigma$ to be a limit linear series implies that either $h^0(X, \eta_X + 2x - q) = 2$ or $h^0(X, \eta_X + p- q) = 2$ if $p\in X$. As $(X, x)\in \CMcal{M}_{g-1,1}$ was chosen generic there is no such $q\in X$.  
\end{proof}

\begin{cor}
  We have
  \begin{align*}
    & a_0^w = \frac{g^2 + 3g - 2}{8},\qquad b_0^w = g-1,\\
    & a_0^c = \frac{9g^2 + 59g - 50}{8},\qquad b_0^c = 24g-22.
  \end{align*}
  \label{cor:alpha0_and_beta0}
\end{cor}

\begin{proof}
  Intersecting both sides of~\eqref{eq:notation_for_coefficietns} with $F_0$ and $G_0$ and using Lemma~\ref{lemma:Cau_and_W_cdot_F0_and_G0} we find that
  \begin{align*}
    l^w -12a_0^w + a_1^w = 0,\\
    3l^w - 12a_0^w -  12b_0^w + 3a_{g-1}^w = 0
  \end{align*}
  Substituting values of $l^w$ and $a_1^w, a_{g-1}^w$ computed in Proposition~\ref{prop:relation_between_W_and_lambda} and Proposition~\ref{prop:alphaj_for_j>0} we obtain two relations for $a_0^w$ and $b_0^w$ that define values of these coefficients.
  
  The case of $a_0^c$ and $b_0^c$ can be handled using Lemma~\ref{lemma:Cau_and_W_cdot_F0_and_G0}, Proposition~\ref{prop:relation_between_Cau_and_lambda} and Proposition~\ref{prop:alphaj_for_j>0} in the same way.
\end{proof}

\section{Test relation for coefficients at $\lambda$ in Theorem~\ref{thm:caustic_and_base_point_divisors_formula}}\label{sec:verification_section}

We devote the last section of our work to present a relation between $[\Cau]$, $[\W]$ and $\lambda$ in $\Pic(\So)$ obtained by pulling back a relation from the moduli space of quadratic differentials. We will use the notation of Section~\ref{sec:construction_over_the_smooth_part}, in particular we will use the notion of $\tSop$ and the vector bundle $E\to \tSop$. Let
\begin{align*}
  & \Wop := \{ (C, \eta, p)\in \tSop\ \mid\ \exists q\in C:\ h^0(C, \eta + p - q) = 2 \},\\
  & \Caup := \text{closure of }\{ (C, \eta, p)\in \tSop\ \mid\ \exists q\in C:\ \eta + p\geq 3q,\ p\neq q \}.
\end{align*}
Recall that $l: \tSop\to \So$ is the forgetful map. We have $l_*\Wop = 2\W$ and $l_*\Caup = \Cau$.

Let $\pi_1: \Mp\times_{\M}\Mp\to \Mp$ be the projection onto the first factor and $\Delta\subset  \Mp\times_{\M}\Mp$ be the diagonal. Let $\Q_g$ be the total space of rank $3g-2$ vector bundle $(\pi_1)_* \omega_{\pi_1}^{\otimes 2}(\Delta)$ on $\Mp$. The space $\Q_g$ can be thought as the moduli space of pairs $(C, p, \omega)$ where $(C,p)$ is a marked curve of genus $g$ and $\omega$ is a quadratic differential with a simple pole at $p$. Introduce the divisor $D_{\deg}\subset \mathbb P\Q_g$ defined by
\begin{equation}
  D_{\deg} := \{(C,p,\omega)\in \mathbb P\Q_g\ \mid\ \omega\text{ has a zero of order bigger then $1$}\}.
  \label{eq:def_of_Ddeg}
\end{equation}
Let $\psi_1^q\in \Pic(\mathbb P \Q_g)\otimes \mathbb Q$ be the pullback of the tautological class from $\Mp$ and let $\L\to \mathbb P\Q_g$ be the tautological bundle.

\begin{prop}
  The following formula holds in $\Pic(\mathbb P\Q_g)\otimes \mathbb Q$:
  \begin{equation*}
    D_{\deg} \equiv 72\lambda - (10g - 8)c_1(\L) + 4\psi_1^q,
  \end{equation*}
  where $\lambda$ is the pullback of the Hodge class from $\M$.
  \label{prop:Ddeg_quadratic_diff}
\end{prop}

\begin{proof}[Sketch of the proof]
  Let $\Mpp := \Mp\times_{\Mp}\Mp$ and $\C = \mathbb P \Q_g\times_{\Mp}\Mpp$ be the universal curve over $\mathbb P\Q_g$. Let $\pr_1,\pr_2$ denote projections onto the first and the second factors of $\mathbb P \Q_g\times_{\Mp}\Mpp$. Let $\widetilde{\mathbb P\Q}_g := \div(\pr_1^*\L\to \pr_2^*\omega_{\pi_1}^{\otimes 2}(\Delta))$. Alternatively, we can write
  \begin{equation*}
    \widetilde{\mathbb P\Q}_g = \{(C,p,\omega, q)\ \mid\ (C,p,\omega)\in \mathbb P\Q_g,\ \omega(q) = 0\}.
  \end{equation*}
  Then the projection $\widetilde{\mathbb P\Q}_g\to \mathbb P\Q_g$ has a simple ramification over $D_{\deg}$, thus we can write applying the adjunction formula:
  \begin{equation*}
    D_{\deg} \equiv (\pr_1)_*\left( \widetilde{\mathbb P\Q}_g\cdot (\widetilde{\mathbb P\Q}_g + \omega_{\pr_1}) \right) \equiv  (\pr_1)_*\left( \widetilde{\mathbb P\Q}_g\cdot (\widetilde{\mathbb P\Q}_g + \pr_2^*(\omega_{\pi_1})) \right).
  \end{equation*}
  Note that this formula holds only if $\widetilde{\mathbb P\Q}_g$ is smooth; we omit a proof of this fact. Expanding the divisor class of $\widetilde{\mathbb P\Q}_g$ as $\pr_2^*(2\omega_{\pi_1}+\Delta) - \pr_1^*c_1(\L)$ and using Mumford's formula for the first kappa class we obtain
  \begin{align*}
    D_{\deg} & \equiv (\pr_1)_*\Bigl( ( \pr_2^*(2\omega_{\pi_1} + \Delta) - \pr_1^*c_1(\L) \Bigr)\cdot \Bigl ( \pr_2^*(3\omega_{\pi_1} + \Delta) - \pr_1^*c_1(\L) )  \Bigr)\\
    & \equiv  72\lambda - (10g - 8)c_1(\L) + 4\psi_1^q.
  \end{align*}
\end{proof}

Let us show how to use Proposition~\ref{prop:Ddeg_quadratic_diff} to establish a relation between $[\Cau]$, $[\W]$ and $\lambda$. Let $(C,\eta,p)\in \tSop$ and let $\sigma_0,\sigma_{\infty}\in E|_{(C,\eta,p)}$ be two elements that span the fiber $E|_{(C,\eta,p)}$. Let $z: \mathbb P H^0(C,\eta + p)\to \mathbb P^1$ be the meromorphic coordinate such that $z(\sigma_0) = 0, z(\sigma_{\infty}) = \infty, z(\sigma_0 + \sigma_{\infty}) = 1$. Then the composition map $C\stackrel{|\eta + p|}{\longrightarrow} H^0(C,\eta + p) \stackrel{z}{\longrightarrow} \mathbb P^1$ is a meromorphic function that we denote simply by $\frac{\sigma_0}{\sigma_{\infty}}$. Recall that $E|_{(C,\eta,p)} \simeq H^0(C,\eta + p)$. Using this isomorphism we identify $\sigma_{\infty}^2$ with a meromorphic differential. Then $\omega_{\sigma_0,\sigma_{\infty}} := d\left( \frac{\sigma_0}{\sigma_{\infty}} \right)\cdot \sigma_{\infty}^2$ is a quadratic differential with at most simple pole at $p$ and no other poles. If we choose another basis $\sigma_0,\sigma_1$ then the differential $\omega_{\sigma_0,\sigma_{\infty}}$ will be multiplied by the determinant of the matrix of the change between bases. Thus the map $\xi: \tSop\to \So\times_{\M}\mathbb P\Q_g$ that sends $(C,\eta,p)$ to $(C,\eta,p, \omega_{\sigma_0,\sigma_{\infty}})$ is defined independently of $\sigma_0,\sigma_{\infty}$. Moreover, the behavior of $\omega_{\sigma_0,\sigma_{\infty}}$ under a change of the basis and~\eqref{eq:first_chern_class_of_spin_bundle} imply that
\begin{equation}
  \begin{split}
    & \xi^*c_1(\L) = c_1(E) - \frac{\lambda}{2},\\
    & \xi^* \psi_1^q = \psi_1,
  \end{split}
  \label{eq:pullback_under_xi}
\end{equation}
where $\psi_1$ was defined in~\eqref{eq:definition_of_psi1}. We claim that
\begin{equation}
  \xi^*D_{\deg} = \Wop + \Caup.
  \label{eq:pullback_of_Ddeg}
\end{equation}
To handle this claim we first note that $\div(\omega_{\sigma_0,\sigma_{\infty}}) + 2p$ is exactly the ramification divisor of the map $C\stackrel{|\eta + p|}{\longrightarrow} \mathbb PH^0(C,\eta + p)$. It follows that $\xi^{-1}(D_{\deg})$ coincides with the locus where two ramification points coalesce, thus $\xi^{-1}(D_{\deg}) = \Wop\cup\Caup$. The fact that $\Wop$ and $\Caup$ appears in the pullback of $D_{\deg}$ with multiplicity $1$ can be deduced from Corollary~\ref{cor:cor_tSopt_is_smooth}.

Combining~\eqref{eq:pullback_under_xi},~\eqref{eq:pullback_of_Ddeg} and Proposition~\ref{prop:Ddeg_quadratic_diff} we obtain  
\begin{equation*}
  \Wop + \Caup \equiv (5g+68)\lambda - (10g - 8)c_1(E) + 4\psi_1.
\end{equation*}
Using~\eqref{eq:first_chern_class_of_E} and pushing this relation forward to $\So$ we find that
\begin{equation}
  2\W + \Cau \equiv (5g^2 + 95g - 70) \,\lambda.
  \label{eq:verifying_end}
\end{equation}
This relation is consistent with Theorem~\ref{thm:caustic_and_base_point_divisors_formula}.


\begin{thebibliography}{00}




\bibitem{ATJ}
M. Atiyah, {\it Riemann surfaces and spin structures,} Ann. Scient. Ec. Norm. Sup.
4, 47--62 (1971).

\bibitem{BOSSY}
  Corentin Boissy, {\it Connected components of the strata of the moduli space of meromorphic differentials,} (2017).



\bibitem{COR}
M. Cornalba, {\it Moduli of curves and theta-characterstics,} Lectures on Riemann surfaces (Trieste), 560--589 (1987).


\bibitem{DK}
  I. Dolgachev and V. Kanev, {\it Polar covariants of plane cubics and quartics}, Advances in Mathematics
98, 216-301, (1993).




\bibitem{LimLin}
  D. Eisenbud and J. Harris, {\it Limit linear series: Basic theory}, Inventiones Math. 85, 337-371, (1986).



\bibitem{KonZor}
A. Eskin, M. Kontsevich, A. Zorich, {\it Sum of Lyapunov exponents of the Hodge bundle with respect to the Teichmüller geodesic flow,} Publications mathématiques de l'IHES,
Volume 120, Issue 1 , 207--333 (2014).


\bibitem{FARo}
G. Farkas, A. Verra, {\it The geometry of the moduli space of odd spin curves,} Annals of Mathematics 180 (2014), 927-970.



\bibitem{ConComp}
M. Kontsevich, A. Zorich, {\it Connected components  of  the moduli spaces of Abelian differentials 
with prescribed singularities, Inventiones Mathematicae,} 153, 631--678 (2003).



\bibitem{KZ}
D. Korotkin, P. Zograf, {\it Tau function and moduli of differentials,} Math. Res. Lett. 18, no.3,
447--458 (2011).




\bibitem{MUM}
D. Mumford, {\it Theta-characteristics of an algebraic curve,} Ann. Scient. Ec.
Norm. Sup. 2, 181--191 (1971).






\end{thebibliography}
\end{document}